\documentclass[final,leqno]{siamltex704}
\usepackage{amsmath,amsfonts,amssymb}
\newtheorem{teor}{Theorem}[section]

\newtheorem{prop}[teor]{Proposition}

\newtheorem{defi}[teor]{Definition}
\newtheorem{nota}[teor]{ {\it Remark}}
\numberwithin{equation}{section}
\newcommand{\R}{\mathbb{R}}

\newcommand{\N}{\mathbb{N}}

\newcommand{\Om}{\Omega}
\newcommand{\w}{\omega}
\newcommand{\G}{\Gamma}
\newcommand{\ep}{\varepsilon}

\newcommand{\di}{\textsf{d}}
\newcommand{\wt}{\widetilde}

\DeclareMathOperator{\car}{card}
\DeclareMathOperator{\cls}{cls}
\newcommand{\n}[1]{\| #1 \|}

\title{Exponential ordering for nonautonomous \\ neutral functional differential
equations}

\author{Sylvia Novo, Rafael Obaya\thanks{Dpto. de Matem\'{a}tica Aplicada, E.T.S. de Ingenieros Industriales, Universidad de
        Va\-lla\-do\-lid, 47011 Valladolid, Spain ({\tt sylnov@wmatem.eis.uva.es}}, {\tt rafoba@wmatem.eis.uva.es}). Partly supported
        by Junta de Castilla y Le\'{o}n  and
        MEC under projects VA024A06 and MTM2008-00700. \and V\'{\i}ctor M. Villarragut\thanks{Dpto. de Matem\'{a}tica Aplicada, E.T.S. de Ingenieros Industriales, Universidad de
        Va\-lla\-do\-lid, 47011 Valladolid, Spain ({\tt vicmun@wmatem.eis.uva.es}}). Partly supported by MEC under project MTM2008-00700 and MICINN under FPU grant AP2006-00874. }

\begin{document}

\vspace*{-1.5cm}

\noindent\begin{tabular}{|l|}
  \hline
  S. Novo, R. Obaya, V. M. Villarragut,\\
  Exponential ordering for nonautonomous neutral functional differential
equations,\\
  SIAM Journal on Mathematical Analysis, Volume 41,
  Number 3, 2009, Pages 1025-1053.\\
  https://doi.org/10.1137/080744682\\
  \copyright ~Society for Industrial and Applied Mathematics\\
  \hline
\end{tabular}

\vspace{1cm}

\maketitle

\begin{abstract}
We study monotone skew-product semiflows generated by families of nonautonomous neutral functional differential equations with infinite delay and stable D-operator, when the exponential ordering is considered. Under adequate hypotheses of  stability for the order on bounded sets, we show that the omega-limit sets are copies of the base to explain the long-term behavior of the trajectories.
The application to the study of the amount of material within the
compartments of a neutral compartmental system with infinite delay, shows the improvement with respect to the standard ordering.
\end{abstract}



\pagestyle{myheadings} \thispagestyle{plain} \markboth{S. NOVO, R. OBAYA AND
V.M. VILLARRAGUT}{EXPONENTIAL ORDERING NFDS}

\section{Introduction}\label{intro}
Dynamical methods for monotone autonomous differential equations have been extensively studied for several decades (see Hirsch~\cite{H1988}, Matano \cite{MA1984}, Pol\'{a}\v{c}ik~\cite{PO1989} and Smith~\cite{smith1995} among many others). A critical problem in this theory concerns the long-term behavior of relatively compact trajectories. In an adequate dynamical scenario these papers prove that the relatively compact trajectories of a strongly monotone semiflow converge  generically to the set of equilibria. This property, referred to as {\em generic quasi-convergence\/}, has a significative theoretical and practical interest.

Following the above theoretical framework, Smith and Thieme~\cite{SMTH1990, SMTH1991} studied the dynamical properties of semiflows induced by functional differential equations with finite delay, which are monotone for the exponential ordering. This ordering is technically more complicated than the
standard one and can be employed in several important situations where the standard quasimonotone condition fails, but a new nonstandard monotone condition is satisfied. Assuming that the vector field satisfies a strong version of the nonstandard monotonicity condition and a supplementary irreducibility assumption, Smith and Thieme proved that the induced semiflow is eventually strongly monotone in the phase space of Lipschitz functions, and strongly order preserving in the phase space of continuous functions, from which they deduced the quasi-convergence property. These results have been extended in several directions: Krisztin and Wu~\cite{KWU1996} studied the dynamical properties of scalar neutral functional differential equations with finite delay which induce a monotone semiflow for the exponential ordering, and Wu and Zhao~\cite{WZ2002} analyzed the same questions in a class of evolutionary equations with applications to reaction-diffusion models.

More recently a new dynamical theory for deterministic and random monotone nonautonomous differential equations has been developed (see Shen and Yi~\cite{shyi}, Jiang and Zhao~\cite{jizh}, Novo, Obaya and Sanz~\cite{NOS2007} and Chueshov~\cite{chue} among others). Now the generic quasi-convergence property fails in general and new dynamical situations have to be considered. Assuming adequate hypotheses including the boundedness, relative compactness and uniform stability of the trajectories, this theory ensures the convergence of the orbits to solutions which reproduce exactly the dynamics exhibited by the time variation of the equation.  The aim of this paper is to extend the applicability of this theory to functional differential equations (FDEs) with infinite delay and  neutral functional differential equations (NFDEs) with infinite delay and $D$-stable operator, which are monotone for the exponential ordering. The obtained results are applied to the study of the long-term behavior of the amount of material within the compartments of some classes of nonautonomous and closed compartmental systems extensively studied in the literature.

Compartmental systems have
been used as mathematical models for the study of the dynamical
behavior of many processes in biological and physical sciences,
which depend on local mass balance conditions (see
Jacquez~\cite{JJ1996}, Jacquez and Simon~\cite{JS1993,JS2002}, and the
references therein). Some initial results for models described by
FDEs with finite and infinite delay can be found in
Gy\"{o}ri~\cite{gyori} and Gy\"{o}ri and Eller~\cite{gyell}. NFDEs represent systems
where the compartments produce or swallow material. The papers by Arino and Haourigui~\cite{ARI1987}, Gy\"{o}ri and
Wu~\cite{gyoriwu}, Wu and Freedman~\cite{WF1991} and Mu\~{n}oz-Villarragut, Novo and Obaya~\cite{MNO2008} are significant predecessors in the use of monotone methods to analyze this kind of problems. On the other hand, supplementary results for the scalar NFDE
\begin{equation}\label{escalar}
\frac{d}{dt} [x(t)-c(t)\,x(t-\tau)]=-h(t,x(t))+h(t-\sigma,x(t-\sigma))
\end{equation}
have been obtained. Arino and Bourad~\cite{AB1990} and Wu~\cite{WU1991} studied this equation for $\tau=\sigma$ and  time independent, time-periodic or almost-periodic $c(t)$, $h(t,x)$.
Krisztin~\cite{KR1992} observed that if $c=\sqrt{2}-1$, $\tau=\pi/4$, $\sigma=7\pi/4$ and $h(t,x)=x$, then equation~\eqref{escalar} has the periodic solution $x(t)=\sin (t)$, and hence the convergence to equilibria is not always true when $\tau\neq \sigma$. Krisztin and Wu~\cite{KWU1996} pointed out that, for some values $\tau\neq \sigma$ with $c(t)$, $h(t,x)$ time-periodic, equation~\eqref{escalar}  generates a monotone semiflow for the exponential ordering, which becomes eventually strongly monotone in the phase space of Lipschitz continuous functions. In this space the asymptotic periodicity of the solutions is deduced.

In this paper we extend the conclusions of~\cite{KWU1996} to nonautonomous closed systems of NFDEs with finite or infinite delay.  The total mass in closed compartmental systems defines an invariant function, which is uniformly continuous for the metric on bounded sets. This implies a property of stability for pairs of trajectories starting at ordered initial data. We will refer to this property as {\em stability for the order on bounded sets\/}; this will be the hypothesis we assume throughout the paper, determining the complete presentation of the theory. It is important to mention that in the case of the exponential ordering, this is a weak version of stability and it is not clear how to deduce the standard stability of any particular trajectory.

We assume some recurrence properties on the temporal variation of the NFDE; thus, its solutions induce a skew-product semiflow with a minimal flow on the base. In particular, the uniform almost-periodic and almost-automorphic cases are included in this formulation. The skew-product formalism permits the analysis of the dynamical properties of the trajectories using methods of ergodic theory and topological dynamics. In this paper we study nonautonomous NFDEs with infinite delay and autonomous $D$-stable operator, satisfying the same hypotheses considered in~\cite{MNO2008}, which in particular imply the invertibility of $D$. We introduce the phase space $BU\subset C((-\infty,0],\R^m)$ of bounded and uniformly continuous functions with the supremum norm where the standard theory provides existence, uniqueness and continuous dependence of the solutions. In our present setting  every bounded trajectory is relatively compact for the compact open topology, the restriction of the semiflow to its omega-limit set is continuous for the metric and it admits a flow extension. We also prove that every solution of the NFDE lying on an omega-limit set is differentiable everywhere. In addition, if two functions belong to the same omega-limit set and they are close in metric then the same happens to their derivatives, due to the continuous variation of the vector field on metric compact sets. From this, for each pair of elements $x$, $y$ of each omega-limit set, we deduce the existence of a supremum and an infimum, which depend continuously on $x$, $y$ for the metric.  These facts are essential in the proofs of the main results of the paper.

We now briefly describe the structure of the paper. Basic notions in topological dynamics, used throughout the rest of the sections, are stated in Section~\ref{prelim}. In Section~\ref{expon} we define the exponential ordering in the space $BU$ to be used for the study of FDEs and NFDEs with infinite delay. For FDEs two cases arise: in Case~I we use the exponential ordering on a compact subinterval of $(-\infty,0]$ and the standard ordering on its complementary subset, and  in Case~II the exponential ordering is used on the complete semi-interval $(-\infty,0]$. Case~III refers to NFDEs, where we consider the exponential ordering on $(-\infty,0]$. The theory provides different dynamical consequences for each choice. We formulate a nonstandard monotone condition for the vector field which implies the monotonicity of the corresponding semiflow.

Section~\ref{omegalimit} deals with the study of the structure of the omega-limit set of a relatively compact trajectory which is uniformly stable for the order in bounded sets. We make use of the arguments included in Section 3 of Novo, Obaya and Sanz~\cite{NOS2007}; here some supplementary work is required because the conditions of stability are weaker. We prove that this compact set is uniformly stable for the order in bounded sets and the restriction of the semiflow to it is uniformly stable. We finally prove that this compact set is a minimal set assuming that the trajectory starts at a Lipschitz initial value, as a supplementary condition in Cases~II and~III.

Section~\ref{copiadelabase} contains the main statement of the paper, which ensures the convergence of the trajectories and follows the arguments of Jiang and Zhao~\cite{jizh} and Novo, Obaya and Sanz~\cite{NOS2007}. We describe an adequate dynamical scenario where the trajectories with  Lipschitz initial data are relatively compact and uniformly stable for the order on bounded sets. We assume that the vector field satisfies a strong version of the nonstandard monotonicity condition for each component. We now fix a relatively compact trajectory uniformly stable for the order in bounded sets and we prove that its omega-limit set is a minimal set given by a 1-cover of the base flow. Again, in Cases II and III the Lipschitz character of the initial data is assumed. Notice that no irreducibility condition  is required and hence the conclusions can be applied to more general problems under natural physical conditions.

In Section~\ref{example} we apply the previous results to some nonautonomous and closed compartmental systems with finite or infinite delay. We extend a part of the dynamical description detailed in~\cite{MNO2008} to some models which are non-monotone for the standard ordering. In particular, we provide a nonautonomous vectorial version of the results of Kristzin and Wu~\cite{KWU1996}.
\section{Some preliminaries}\label{prelim}
Let $(\Om,d)$ be a compact metric space. A real {\em continuous
flow \/}  $(\Om,\sigma,\R)$ is defined by a continuous mapping
$\sigma: \R\times \Om \to  \Om,\; (t,\w)\mapsto \sigma(t,\w)$
satisfying

\begin{enumerate}
\renewcommand{\labelenumi}{(\roman{enumi})}
\item $\sigma_0=\text{Id},$
\item $\sigma_{t+s}=\sigma_t\circ\sigma_s$ for each $s$, $t\in\R$,
\end{enumerate}
where $\sigma_t(\w)=\sigma(t,\w)$ for all $\w \in \Om$ and $t\in
\R$. The set $\{ \sigma_t(\w) \mid t\in\R\}$ is called the {\em
orbit\/} or the {\em trajectory\/} of the point $\w$. We say that
a subset $\Om_1\subset \Om$ is {\em $\sigma$-invariant\/} if
$\sigma_t(\Om_1)=\Om_1$ for every $t\in\R$.  A subset
$\Om_1\subset  \Om$  is called {\em minimal \/} if it is compact,
$\sigma$-invariant and its only nonempty compact
$\sigma$-invariant subset is itself. Every compact and
$\sigma$-invariant  set contains a minimal subset; in particular
it is  easy to prove that a compact $\sigma$-invariant subset is
minimal if and only if every trajectory is dense. We say that the
continuous flow $(\Om,\sigma,\R)$ is {\em recurrent\/} or {\em
minimal\/} if $\Om$ is minimal.

The flow $(\Om,\sigma,\R)$ is {\em distal\/} if for any two
distinct points $\w_1,\,\w_2\in\Om$ the orbits keep at a positive
distance, that is, $\inf_{t\in
\R}d(\sigma(t,\w_1),\sigma(t,\w_2))>0$. The flow $(\Om,\sigma,\R)$
is {\em almost periodic\/} when for every $\varepsilon
> 0 $ there is a $\delta >0$ such that, if $\w_1$, $\w_2\in\Om$
with $d(\w_1,\w_2)<\delta$, then
$d(\sigma(t,\w_1),\sigma(t,\w_2))<\varepsilon$ for every $t\in
\R$. If $(\Om,\sigma,\R)$ is almost periodic, it is distal. The
converse is not true; even if $(\Om,\sigma,\R)$ is minimal and
distal, it does not need to be almost periodic. For the basic
properties of almost periodic and distal flows we refer the reader
to Ellis~\cite{elli} and Sacker and Sell~\cite{sase1}.

A {\em flow homomorphism\/} from another continuous flow
$(Y,\Psi,\R)$ to $(\Om,\sigma,\R)$ is a continuous map $\pi\colon
Y\to \Om$ such that $\pi(\Psi(t,y))=\sigma(t,\pi(y))$ for every
$y\in Y$ and $t\in\R$. If $\pi$ is also bijective, it is called a
{\em flow isomorphism\/}. Let $\pi:Y \to \Om$ be a surjective flow
homomorphism and suppose $(Y,\Psi,\R)$ is minimal (then, so is
$(\Om,\sigma,\R)$). $(Y,\Psi,\R)$ is said to be an {\em almost
automorphic extension\/} of $(\Om,\sigma,\R)$ if there is $\w\in
\Om$ such that $\car (\pi^{-1}(\w))=1$. Then, actually $\car
(\pi^{-1}(\w))=1$ for $\w$ in a residual subset $\Om_0\subseteq
\Om$; in the nontrivial case $\Om_0\subsetneq \Om$ the dynamics
can be very complicated. A minimal flow $(Y,\Psi,\R)$ is {\em
almost automorphic\/} if it is an almost automorphic extension of
an almost periodic minimal flow $(\Om,\sigma,\R)$. We refer the
reader to the work of Shen and Yi~\cite{shyi} for a survey of
almost periodic and almost automorphic dynamics.

Let $E$ be a complete metric space and  $\R^+=\{t\in\R\,|\,t\geq
0\}$. A {\em semiflow} $(E,\Phi,\R^+)$ is determined by a
continuous map $\Phi: \R^+\times E \to E,\; (t,x)\mapsto
\Phi(t,x)$ which satisfies

\begin{enumerate}
\renewcommand{\labelenumi}{(\roman{enumi})}
\item $\Phi_0=\text{Id},$
\item $\Phi_{t+s}=\Phi_t \circ \Phi_s\;$ for all  $\; t$, $s\in\R^+,$
\end{enumerate}
where $\Phi_t(x)=\Phi(t,x)$ for each $x \in E$ and $t\in \R^+$.
The set $\{ \Phi_t(x)\mid t\geq 0\}$ is the {\em semiorbit\/} of
the point $x$. A subset  $E_1$ of $E$ is {\em positively
invariant\/} (or just $\Phi$-{\em invariant\/}) if
$\Phi_t(E_1)\subset E_1$ for all $t\geq 0$. A semiflow
$(E,\Phi,\R^+)$ admits a {\em flow extension\/} if there exists a
continuous flow  $(E,\wt \Phi,\R)$ such that $\wt
\Phi(t,x)=\Phi(t,x)$ for all $x\in E$ and $t\in\R^+$. A compact
and positively invariant subset admits a flow extension if the
semiflow restricted to it admits one.

Write $\R^-=\{t\in\R\,|\,t\leq 0\}$. A {\em backward orbit\/} of a
point  $x\in E$ in the semiflow $(E,\Phi,\R^+)$ is a continuous
map  $\psi:\R^-\to E$ such that $\psi(0)=x$ and for each $s\leq 0$
it holds that  $\Phi(t,\psi(s))=\psi(s+t)$ whenever $0\leq t\leq
-s$. If for $x\in E$ the semiorbit $\{\Phi(t,x)\mid t\ge 0\}$ is
relatively compact, we can consider the {\em omega-limit set\/} of
$x$,
\[
\mathcal{O}(x)=\bigcap_{s\ge 0}{\rm closure}{\{\Phi(t+s,x)\mid
t\ge 0\}}\,,
\]
which is a nonempty compact connected and $\Phi$-invariant set.
Namely, it consists of the points $y\in E$ such that $y=\lim_{\,n\to
\infty} \Phi(t_n,x)$ for some sequence $t_n\uparrow \infty$. It is
well-known that every $y\in\mathcal{O}(x)$ admits a backward orbit
inside this set. Actually, a compact positively invariant set $M$
admits a flow extension if every point in $M$ admits a unique
backward orbit which remains inside the set $M$ (see Shen and
Yi~\cite{shyi}, part~II).

A compact positively invariant set $M$ for the semiflow
$(E,\Phi,\R^+)$ is {\em minimal\/} if it does not contain any
other nonempty compact positively invariant set than itself. If
$E$ is minimal, we say that the semiflow is minimal.

A semiflow is {\em of skew-product type\/} when it is defined on a
vector bundle and has a triangular structure; more precisely, a
semiflow $(\Om\times X,\tau,\,\R^+)$ is a {\em skew-product\/}
semiflow over the product space $\Om\times X$, for a compact
metric space $(\Om,d)$ and a complete metric space
$(X,\textsf{d})$, if the continuous map $\tau$ is as follows:
\begin{equation}\label{skewp}
 \begin{array}{cccl}
 \tau \colon  &\R^+\times\Om\times X& \longrightarrow & \Om\times X \\
& (t,\w,x) & \mapsto &(\w{\cdot}t,u(t,\w,x))\,,
\end{array}
\end{equation}
where $(\Om,\sigma,\R)$ is a  real continuous flow
$\sigma:\R\times\Om\rightarrow\Om$, $\,(t,\w)\mapsto \w{\cdot}t$, called
the  {\em base flow\/}. The skew-product semiflow~\eqref{skewp} is
{\em linear\/} if $u(t,\w,x)$ is linear in $x$ for each
$(t,\w)\in\R^+\times\Om$.
\section{Exponential ordering}\label{expon}
 We will study nonautonomous NFDEs with infinite delay and autonomous and
 $D$-stable operator. Firstly, we describe the conditions on the vector field and the neutral operator $D$ assumed throughout the paper.

We consider the Fr\'{e}chet space
$X=C((-\infty,0],\R^m)$ endowed with the compact-open topology,
i.e.~the topology of uniform convergence over compact subsets,
which is a metric space for the distance
\[
\textsf{d}(x,y)=\sum_{n=1}^\infty \frac{1}{2^n}\frac
{\n{x-y}_n}{1+\n{x-y}_n}\,,\quad x,y\in X\,,
\]
where $\n{x}_n=\sup_{s\in[-n,0]}\n{{x(s)}}$, and $\n{\cdot}$ denotes
the maximum norm in $\R^m$.

\vspace{.2cm}
Let $BU\subset X$ be the Banach space
\[BU=\{x\in X\mid  x \text{ is bounded and
  uniformly continuous}\}\]
 with the supremum norm  $\n{x}_\infty=\sup_{s\in(-\infty,0]} \n{x(s)}$.
Given $k>0$,  we will denote
\[B_k=\{x\in BU \mid \n{x}_\infty \leq k\}\,.\] As usual, given
$I=(-\infty,a]\subset\R$,  $t\in I$ and a continuous function
$x:I\to\R^m$, $x_t$ will denote the element of $X$ defined by
$x_t(s)=x(t+s)$ for $s\in (-\infty,0]$.

\vspace{.2cm}
Let $D\colon BU \to \R^m$ be a linear operator satisfying the
hypotheses:

\begin{itemize}
\item[(D1)] $D$ is linear and continuous for the norm.
\item[(D2)] For each $k>0$, $D\colon B_k\to \R^m$ is continuous when we take
the restriction of the compact-open topology to $B_k$, i.e. if
$x_n\stackrel{\textsf{d}\;}\to x$ as $n\to\infty$ with $x_n$,
$x\in B_k$, then $\lim_{\,n\to\infty}Dx_n=Dx$.
\item[(D3)] $D$ is atomic at $0$ (see definition in Hale~\cite{hale} or
Hale and Verduyn Lunel~\cite{hale2}).
\item[(D4)] $D$ is {\em stable}, i.e. there is a
continuous function $c\in C([0,\infty),\R^+)$ with
$\lim_{\,t\to\infty}c(t)=0$ such that, for each $\varphi\in BU$
with $D\varphi=0$, the solution of
\[\left\{
\begin{array}{ll}
Dx_t=0\,, & t\geq 0 \\
x_0=\varphi\,,
\end{array} \right.
\] satisfies $\n{x(t)}\leq c(t)\,\n{\varphi}_\infty$
for each $t\geq 0$.
\end{itemize}

\vspace{.2cm}
From (D1)-(D2), as shown in Mu\~{n}oz-Villarragut, Novo and Obaya~\cite{MNO2008}, if $x\in BU$
\[Dx=\int_{-\infty}^0 [d\mu(s)]\,x(s)\]
where $\mu=[\mu_{ij}]_{i,j\in\{1,\ldots,m\}}$ and $\mu_{ij}$ is a real regular Borel
measure with finite total variation
$|\mu_{ij}|(-\infty,0]<\infty$, for all $i$, $j\in\{1,\ldots,m\}$. From (D3), without
loss of generality, we may assume that
\begin{equation}\label{Dformula} Dx=x(0)-\int_{-\infty}^0
[d\nu(s)]\,x(s)
\end{equation}
where $\nu=[\nu_{ij}]_{i,j\in\{1,\ldots,m\}}$, $\nu_{ij}$ is a
real regular Borel measure with finite total variation, and
$|\nu_{ij}|(\{0\})=0$ for all $i$, $j\in\{1,\ldots,m\}$. For any measurable set $E\subset (-\infty,0]$ we will
denote by $|\nu|(E)$ the $m\times m$ matrix
$[\,|\nu_{ij}|(E)\,]$ and by $\n{\nu}_\infty(E)$ the
corresponding matricial~norm. We define the linear operator
\begin{equation}\label{Dhat}
\begin{array}{lcclcl}
\widehat{D} \colon &BU &\longrightarrow & BU &&\\
& x & \mapsto &\widehat{D}x\colon(-\infty,0]& \to  &\R^m \\
&  & &   \hspace{1.2cm}s & \mapsto & Dx_s\,,
\end{array}
\end{equation}
that is,
$\widehat{D}x(s)=x(s)-\int_{-\infty}^0[d\nu(\theta)]\,x(\theta+s)$
for each $s\in(-\infty,0]$. From (D4) (see Theorem 3.8
of~\cite{MNO2008}) $\widehat{D}$ is invertible, $\widehat{D}^{-1}$
is bounded for the norm and uniformly continuous when we take the
restriction of the compact-open topology to $B_k$, i.e. given
$\varepsilon >0$ there is a $\delta(k)>0$ such that
{\upshape$\textsf{d}(\widehat{D}^{-1}h_1,\widehat{D}^{-1}h_2)<\varepsilon$}
for all $h_1$, $h_2\in B_k$ with
{\upshape$\textsf{d}(h_1,h_2)<\delta(k)$}. Hence the linear
operator $T\colon BU \to \R^m$, $x\mapsto (\widehat{D}^{-1}x)(0)$,
also satisfies (D1)-(D4), and has a representation
\begin{equation}\label{defiT}
Tx=\int_{-\infty}^0 [d\widehat
\mu(s)]\,x(s)
\,,
\end{equation}
where $\widehat\mu=[\widehat\mu_{ij}]$ and $\widehat\mu_{ij}$ is a real regular Borel
measure with finite total variation.

Let $\mathcal{L}$ be the space of linear operators
$\mathcal{L}=\{D\colon BU\to \R^m\mid \text{ (D1)-(D2) hold}\}$,
which is complete for the operator norm. As a consequence of the next result it can be checked that $\mathcal{U}=\{D\in\mathcal{L}\mid \text{ (D3)-(D4) hold}\}$ is
open in $\mathcal{L}$.

\begin{prop}\label{estableabierto}
Let us assume that $D$ satisfies {\rm (D1)-(D4)} and it is given
by~\eqref{Dformula}. Then, there is an $\ep>0$ such that for any
$\nu^*=[\nu^*_{ij}]$ where $\nu^*_{ij}$ is a real regular Borel measure with finite total variation, $|\nu^*_{ij}|(\{0\})=0$ for each $i$, $j\in{1,\ldots,m}$ and  $\n{\nu-\nu^*}_\infty(-\infty,0]<\ep$, the linear
operator $D^*\colon BU \to \R^m$ given by
\[ D^*x=x(0)-\int_{-\infty}^0 [d\nu^*(s)]\,x(s)\]
is stable.
\end{prop}

\begin{proof}
As shown in Theorem 3.9 of~\cite{MNO2008}, it is enough to check
that the linear operator $\widehat{D^*}\colon BU\to BU$, defined
from $D^*$ as in~\eqref{Dhat}, is invertible and $(\widehat
{D^*})^{-1}$ is continuous for the restriction of the compact-open
topology to $B_k$. Since $\widehat D$ is invertible, we define
$X\colon BU\to BU$ by
\[X=\widehat D^{-1}\circ(\widehat D-\widehat{D^*})\,.\]
Hence, taking $\ep <1/\n{\widehat D^{-1}}$, from $\n{\widehat
D-\widehat {D^*}}<\ep$ we deduce that $\n{X}< 1$, ${\rm Id}-X$ is
invertible and consequently $\widehat {D^*}=\widehat D\circ({\rm Id}-X)$
is also invertible. From $(\widehat {D^*})^{-1}=({\rm Id}-X)^{-1}\circ
\widehat D^{-1}=\left(\sum_{n=0}^{\infty}X^n\right)\circ \widehat
D^{-1}$ it is easy to check the continuity for the restriction of the
compact-open topology to $B_k$, and the proof is finished.
\qquad\end{proof}

Let $(\Om,\sigma,\R)$ be a minimal flow over a compact metric
space $(\Om,d)$ and denote $\sigma(t,\w)=\w{\cdot}t$ for all $\w \in\Om$
and $t\in\R$. Let $F\colon\Om\times BU \to\R^m$, $(\w,x)\mapsto
F(\w,x)$ be a function satisfying the following conditions:

\begin{itemize}
\item[(F1)] $F$ is continuous on $\Om\times BU$ and locally
Lipschitz in $x$ for the norm $\n{\cdot}_\infty$. \item[(F2)] For each
$k>0$, $F(\Om\times B_k)$ is a bounded subset of $\R^m$.
\item[(F3)] For each $k>0$, $F\colon\Om\times B_k\to\R^m$ is
continuous when we take the restriction of the compact-open
topology to $B_k$, i.e.~if $\w_n\to\w$ and
$x_n\stackrel{\textsf{d}\;}\to x$ as $n\to\infty$ with $x_n$, $x\in B_k$,
then $\lim_{n\to\infty}F(\w_n,x_n)=F(\w,x)$.
\end{itemize}

\vspace{.1cm}
\noindent We consider the families of nonautonomous FDEs with infinite
delay
\addtocounter{equation}{1}
\begin{equation}\tag*{(\arabic{section}.\arabic{equation})$_\w$}\label{infdelay}
 z'(t)=F(\w{\cdot}t,z_t)\,, \quad t\geq 0\,,\;\w\in\Omega\,,
\end{equation}
and nonautonomous NFDEs with infinite delay and stable
$D$-operator
\addtocounter{equation}{1}
\begin{equation}\tag*{(\arabic{section}.\arabic{equation})$_\w$}\label{Ninfdelay}
 \frac{d}{dt}Dz_t=F(\w{\cdot}t,z_t)\,, \quad t\geq 0\,,\;\w\in\Om\,,
\end{equation}
which obviously includes~\ref{infdelay} when $\nu\equiv
0$.

From hypothesis (F1), the standard theory of
FDEs (resp. NFDEs) with infinite delay, see Hino, Murakami and Naito~\cite{hino} (resp. Wang and Wu~\cite{WW1985} and
Wu~\cite{jwu1991}),  assures that for each $x\in BU$ and each
$\w\in\Om$ the system~\ref{infdelay} (resp. ~\ref{Ninfdelay})
locally admits a unique solution $z(t,\w,x)$ with initial value
$x$, i.e.~$z(s,\w,x)=x(s)$ for each $s\in (-\infty,0]$. Therefore,
the family~\ref{infdelay} (resp. ~\ref{Ninfdelay}) induces a local
skew-product semiflow
\begin{equation}\label{Ndelaskewcoop}
 \begin{array}{cccl}
 \tau &: \R^+\times\Om\times BU& \longrightarrow & \Om\times BU\\
& (t,\w,x) & \mapsto &(\w{\cdot}t,u(t,\w,x))\,,
\end{array}
\end{equation}
where $u(t,\w,x)\in BU$ and $u(t,\w,x)(s)=z(t+s,\w,x)$ for $s\in
(-\infty,0]$.

As shown in Mu\~{n}oz-Villarragut, Novo and Obaya~\cite{MNO2008}, the change of variable $y=\widehat{D}z$ takes~\ref{Ninfdelay} to
\addtocounter{equation}{1}
\begin{equation}\tag*{(\arabic{section}.\arabic{equation})$_\w$}\label{infdelayt}
 y'(t)=G(\w{\cdot}t,y_t)\,, \quad t\geq 0\,,\;\w\in\Om\,,
\end{equation}
with $G\colon\Om\times BU \to\R^m$, $(\w,x)\mapsto
G(\w,x)=F(\w,\widehat{D}^{-1}x)$. This family induces a local
skew-product semiflow
\begin{equation*}
 \begin{array}{cccl}
 \widehat\tau &: \R^+\times\Om\times BU& \longrightarrow & \Om\times BU\\
& (t,\w,x) & \mapsto &(\w{\cdot}t,\widehat u(t,\w,x))\,,
\end{array}
\end{equation*}
where $\widehat u(t,\w,x)\in BU$ and $\widehat
u(t,\w,x)(s)=y(t+s,\w,x)$ for $s\in (-\infty,0]$, and it is
related to the previous one~\eqref{Ndelaskewcoop} by
\begin{equation}\label{usombrero}
\widehat u(t,\w,x)= \widehat{D}\,u(t,\w,\widehat{D}^{-1}\,x)\,.
\end{equation}

The above change of variable is often used to deduce the regularity of the solutions of~\ref{Ninfdelay}. In particular, if the initial value $x$ is Lipschitz, then it is easy to check that $\widehat{D}x$ is Lipschitz, and as a consequence, $\widehat u(t,\w,\widehat D x)$ and $u(t,\w,x)$ are also Lipschitz.
Besides, as it is well-known~\cite{hale, hale2}, the solution of a NFDE is not necessarily differentiable at every point. Next we
show that this indeed holds provided that the trajectory is bounded
and admits a backward orbit.

\begin{prop}\label{clasec1}
We consider the skew-product semiflow~\eqref{Ndelaskewcoop}
induced by equation~{\rm\ref{Ninfdelay}}. Assume that
$(\w,x)\in\Om\times BU$ admits a backward orbit extension and that
there is a $k_1>0$ such that $u(t,\w,x)\in B_{k_1}$ for each
$t\in\R$. Then $z(t)=z(t,\w,x)$, the solution
of~{\rm\ref{Ninfdelay}} with initial value $x$, belongs to
$C^1(\R,\R^m)$.
\end{prop}

\begin{proof}
From~\eqref{defiT}, we have for each $h\neq 0$
\[\frac{z(t+h)-z(t)}{h}=\int_{-\infty}^0 [ d\widehat
\mu (\theta)]\,\frac{y(t+\theta+h)-y(t+\theta)}{h}\,,\] where $y(t)=y(t,\w,{\widehat D}x)$ is the solution
of~\ref{infdelayt} with initial value ${\widehat D}x$, which is
defined for each $t\in\R$, it belongs to $C^1(\R,\R^m)$ and is bounded by
$k_2=\n{\widehat D}\,k_1$. From (F2) there is a $c>0$ such that  $\n{G(\Om\times
B_{k_2})}\leq c$; then an application of the dominated convergence
theorem shows that $z$ is differentiable,
\begin{equation}\label{zprima}
z'(t)=\int_{-\infty}^0 [d\widehat \mu(\theta)]\,y'(t+\theta)\,,
\end{equation}
and $\n{z'(t)}\leq c\,\n{\widehat \mu}_\infty(-\infty,0]$ for each
$t\in\R$.

Finally, we check the continuity of  $z'$. Given $t_0$ and
$\varepsilon>0$ we find $k>0$ such that $\n{\widehat\mu}_\infty(-\infty,-k]<\varepsilon/(4\,c)$. Since $y'$ is uniformly
continuous on $I_0=[t_0-k-1,t_0+1]$, there is $0<\delta< 1$ such
that $\n{y'(s)-y'(\theta)}<\varepsilon/(2\, \n{\widehat \mu}_\infty[-k,0]+1)$
for each $s$, $\theta\in I_0$ with $|s-u|<\delta$. Hence, from~\eqref{zprima}, we deduce
that if $|t-t_0|<\delta$
\begin{align*}
\n{z'(t)-z'(t_0)}  & \leq  \frac{\varepsilon}{2}+  \left\|\int_{-k}^{0}
\,[d\widehat\mu(\theta)]\,(y'(t+\theta)-y'(t_0+\theta))\,\right\|\\
&  \leq  \frac{\varepsilon}{2}+ \n{\widehat{\mu}}_\infty[-k,0] \sup_{\theta\in[-k,0]}\n{y'(t+\theta)-y'(t_0+\theta)} \leq
 \frac{\varepsilon}{2}+ \frac{\varepsilon}{2}=\varepsilon
\end{align*}
because $t+\theta$ and $t_0+\theta$ belong to $I_0$ for $\theta\in[-k,0]$,
and the proof is finished.
\qquad\end{proof}

Next, we introduce two exponential orderings in the space $BU$ and we will formulate
a simultaneous quasimonotone condition for FDEs and NFDEs.

As usual, in $\R^m$ we take the partial order relation
\[
\begin{split}
 v\le w \quad&\Longleftrightarrow\quad v_j\leq w_j\quad \text{for}\;j=1,\ldots,m\,,\\
 v<w  \quad  &\Longleftrightarrow\quad\,v\,\leq\, w\quad\text{and}
 \quad v_j<w_j\quad\text{for some}\;j\in\{1,\ldots,m\}\,,
\end{split}
\]
and if $x\in BU$ we denote $x\geq 0$ when $x(s)\geq 0$ for
each $s\in(-\infty,0]$.

We write $A\leq B$ for $m\times m$ matrices $A=[a_{ij}]$ and $B=[b_{ij}]$ if and only if  $a_{ij}\leq b_{ij}$ for all $i,j$. Let $A$ be an $m\times m$ quasipositive matrix, that is, there is a $\lambda>0$ such that $A+\lambda I\geq 0$.

\begin{nota}\label{quasipositive}
\rm If $A$ is a quasipositive matrix then $e^{At}$ is a nonnegative matrix with the entries in the main diagonal being strictly positive for each $t\geq 0$.
\end{nota}

We introduce
the new families of positive cones with empty interior on $BU$
\begin{align*} BU^+_{A,r}&=\{
x\in BU\mid x\geq 0 \;\text{ and }\; x(t)\geq e^{A(t-s)}x(s)\quad \text{for }
-r\leq s\leq t\leq 0\}\,,\\
BU^+_{A,\infty}&=\{ x\in BU\mid x\geq 0 \;\text{ and }\;
x(t)\geq e^{A(t-s)}x(s)\quad \text{for }
-\infty<s\leq t\leq 0\}\,,
\end{align*}
which respectively induce the following partial order relations on
$BU$
\[
\begin{split}
 x\le_{A,r} y \;&\Longleftrightarrow\; x\le y
\;\text{ and }\; y(t)-x(t)\geq e^{A(t-s)}(y(s)-x(s))\,,
-r\leq s\leq t\leq 0\,,\\
 x<_{A,r} y \;&\Longleftrightarrow\; x \leq_{A,r} y \quad\text{and}\quad x\neq
 y\,,\;\text{and} \\[.1cm]
 x\le_{A,\infty} y\;&\Longleftrightarrow\; x\le y \;\text{ and }\;
y(t)-x(t)\geq e^{A(t-s)}(y(s)-x(s))\,,
-\infty<s\leq t\leq 0\,,\\
 x<_{A,\infty} y  \;&\Longleftrightarrow\; x \leq_{A,\infty} y\;\;\text{and}\quad x\neq
 y\,.
\end{split}
\]
As in~\cite{SMTH1991} for the case of finite delay, the following result  can be checked.

\begin{nota}\label{lipschitz}
\rm{
A smooth function (resp. a Lipschitz continuous function) $x$ belongs to $BU_{A,r}^+$ if and only if
\[x\geq 0 \quad\text{and}\quad x'(s)\geq A\,x(s)\quad\text{for each (resp. a.e.) }\; s\in[-r,0]\,,\]
and belongs  to $BU_{A,\infty}^+$ if and only if
\[\;x\geq 0 \quad\text{and}\quad x'(s)\geq A\,x(s)\quad\text{for each (resp. a.e.) }\; s\in(-\infty,0]\,.\]}
\end{nota}
\indent In the rest of the paper $\leq_A$ will denote any of the above
order relations and we assume the following quasimonotone
condition

\begin{itemize}
\item[(F4)] If
$x,y\in BU$ with $x\leq_A y$ then
$F(\w,y)-F(\w,x)\geq A\,(Dy-Dx)$ for each $\w\in\Om$\,,
\end{itemize}
to generate a monotone skew-product semiflow in the following cases:

\vspace{.2cm}\begin{tabular}{ll}
        Case I: & FDEs~\ref{infdelay} under (F1)-(F4) and order  $\leq_{A,r}$;\\
        Case II: & FDEs~\ref{infdelay} under (F1)-(F4) and order  $\leq_{A,\infty}$;\\
        Case III: & NFDEs~\ref{Ninfdelay}, under (F1)-(F4), (D1)-(D5) and order  $\leq_{A,\infty}$;\vspace{.2cm}
\end{tabular}

\noindent with $D$ of the form~\eqref{Dformula}, $\nu\equiv 0$ in Cases I-II, and

\begin{itemize}
\item[(D5)] The measures $\nu_{ij}$ in~\eqref{Dformula} are positive for all $i$, $j\in\{1,\ldots,m\}$.
\end{itemize}

\begin{teor}\label{monotone} The skew-product
semiflow~{\rm \eqref{Ndelaskewcoop}} induced
by~{\rm\ref{Ninfdelay}} is monotone, that is,  for each
$\w\in\Omega$ and $x, y\in BU$ such that $x\leq_A y$ it holds
that
\[u(t,\w,x)\leq_A u(t,\w,y)\]
whenever they are defined.
\end{teor}

\begin{proof}
 We omit the proof of Cases I-II which is completely analogous to the one given in Proposition 3.1 of~\cite{SMTH1991}.

Case III. First, we show the result when $D$ has the form
\begin{equation}\label{Dmenosrho}
Dx=x(0)-\int_{-\infty}^{-\rho}[d\nu(s)]\, x(s)
\end{equation}
for some $\rho>0$. Let ${\bf 1}=(1,\ldots,1)^T\in \R^m$. For $\varepsilon>0$ let $z_\varepsilon(t,\w,y)$ be the
solution~of
\[
 \frac{d}{dt}Dz_t=F(\w{\cdot}t,z_t)+\varepsilon \,{\bf 1}\,, \quad t\geq 0\,,\;\w\in\Om\,,
\]
with initial value $y$, i.e. $z_\varepsilon(s,\w,y)=y(s)$ for each
$s\in(-\infty,0]$. There exist $\varepsilon_0>0$ and $T>0$ such
that for each $\varepsilon\in[0,\varepsilon_0)$,
$z(t)=z(t,\w,x)$ and
$y^\varepsilon(t)=z_\varepsilon(t,\w,y)$ are defined on $[0,T]$
and denote
$z^\varepsilon(t)=y^\varepsilon(t)-z(t)$.

Let $t_1\in[0,T]$ be the greatest  time such that
$z^\varepsilon_t\geq_A 0$ (i.e. $z_t\leq_A
y^\varepsilon_t$), that is, $z^\varepsilon(t)\geq 0$ and $z^\varepsilon(t)\geq e^{A(t-s)} z^\ep(s)$ for each $s\leq t\leq t_1$.  We claim that $t_1=T$. Assume on the contrary
that $t_1<T$.  From (F4)
\[  \frac{d}{dt}Dz^\varepsilon_t |_{t=t_1}-ADz_{t_1}^\ep=
[F(\w,y^\varepsilon_{t_1})-F(\w,z_{t_1})-A
D(y^\varepsilon_{t_1}-z_{t_1})]+
\varepsilon\,{\bf 1}\geq \ep\,{\bf 1}\,,\] and there is
 $h>0$ with $h<\rho$ such that
\begin{equation*}
 \frac{d}{dt}
Dz^\varepsilon_t -ADz_t^\ep\geq 0\,,\quad t\in[t_1,t_1+h]\,.
\end{equation*}
Hence, $Dz^\varepsilon_t\geq e^{A(t-s)}Dz^\varepsilon_s$
for each $t_1\leq s\leq t\leq t_1+h$,
and from~\eqref{Dmenosrho}
\[  z^\varepsilon (t)- e^{A(t-s)}
z^\varepsilon(s)  \geq \int_{-\infty}^ {-\rho}
[d\nu(\theta)](z^\varepsilon(t+\theta)-e^{A(t-s)}z^\varepsilon(s+\theta))\,.\]
 Since $h<\rho$ and
$t\in[t_1,t_1+h]$ we have $s+\theta\leq t+\theta\leq t_1+h-\rho\leq t_1$ and
thus $z^\varepsilon(t+\theta)-e^{A
(t-s)}z^\varepsilon(s+\theta)\geq 0$ if $\theta\in(-\infty,-\rho)$.
Consequently, from (D5)
\begin{equation*}
z^\varepsilon (t)- e^{A(t-s)} z^\varepsilon(s)\geq 0
\end{equation*}
if $t_1\leq s\leq t\leq t_1+h$,
which contradicts the definition of $t_1$. Then $t_1=T$ and
letting $\varepsilon$ go to $0$ the result is
proved.

Next, since  $D$ is stable,  from
Proposition~\ref{estableabierto} there is an $\ep_0>0$ such that
\[D_\ep: BU\to\R^m\,,\quad x\mapsto D_\ep x=x(0)-\int_{-\infty}^{-\ep}\,[d\nu(s)]\,x(s)\] is stable for
each $\ep\leq \ep_0$. Then the skew-product semiflow
\begin{equation*}
 \begin{array}{cccl}
 \tau_\ep &: \R^+\times\Om\times BU & \longrightarrow & \Om\times BU\\
& (t,\w,x) & \mapsto &(\w{\cdot}t,u_\ep(t,\w,x))\,,
\end{array}
\end{equation*}
induced by
\addtocounter{equation}{1}
\begin{equation}\tag*{(\arabic{section}.\arabic{equation})$_\w$}\label{Nepsilom}
 \frac{d}{dt}D_\ep z_t=F(\w{\cdot}t,z_t)-A(Dz_t-D_\ep z_t)\,, \quad t\geq 0\,,\;\w\in\Om\,,
\end{equation}
satisfies the corresponding hypotheses {\rm(F1)-(F4)}  for $D_\ep$, and we deduce~that
\[u_\ep(t,\w,x)\leq_A u_\ep (t,\w,y)\]
provided that $x\leq_A y$. The continuous dependence of
solutions of NFDEs, or the continuous dependence of solutions of
FDEs after taking~\ref{Nepsilom} to a FDEs system as above, shows that $\lim_{\,\ep\to 0}
u_\ep(t,\w,x)=u(t,\w,x)$, and we conclude that $u(t,\w,x)\leq_A
u(t,\w,y)$, as stated.
\qquad\end{proof}

Notice that the restriction  provided by Krisztin and Wu in~\cite{KWU1996} on $A$ ($-\mu$ in their paper) to show the theorem in the scalar neutral case with finite delay, has been removed in the previous theorem.
\section{Stability of omega-limit sets}\label{omegalimit}
Throughout this section we will study the omega-limit sets for the
monotone skew-product semiflows~\eqref{Ndelaskewcoop} induced by~\ref{Ninfdelay} in the cases

\vspace{.2cm}\begin{tabular}{ll}
        \!\!Case I: & FDEs~\ref{infdelay} under (F1)-(F4) and order  $\leq_{A,r}$;\\
         \!\!Case II: & FDEs~\ref{infdelay} under (F1)-(F4), (A1) and order  $\leq_{A,\infty}$;\\
         \!\!Case III: & NFDEs~\ref{Ninfdelay}, under (F1)-(F4), (A1), (D1)-(D5) and order  $\leq_{A,\infty}$;\vspace{.2cm}
\end{tabular}

\noindent with $D$ of the form~\eqref{Dformula}, $\nu\equiv 0$ in Cases I-II, and

\vspace{.1cm}
\begin{itemize}
\item[(A1)] All the eigenvalues of $A$ have strictly negative real part.
\end{itemize}
We introduce different concepts of stability for the metric topology of trajectories and compact invariant sets. In Section~\ref{example}, we will illustrate the practical interest of these concepts for the study of compartmental systems which are monotone for the exponential ordering.

\begin{defi}
{\rm Let $K$ be a compact positively $\tau$-invariant set of $\Om\times
BU$. It is said that $(K,\tau,\R^+)$ is  {\em uniformly stable\/}
if for any $\varepsilon>0$ there exists a $\delta>0$, called the
{\em modulus of uniform stability\/}, such that, if $(\w,x),
(\w,y)\in K$ are such that $\di(x,y)< \delta$, then
$\di(u(t,\w,x),u(t,\w,y))\leq \varepsilon$ for all $t\ge 0$.}
\end{defi}

\begin{defi}
{\rm Given $k>0$, a forward orbit $\{\tau(t,\w_0,x_0)\,|\;t\geq 0\}$ of
the skew-product semiflow~\eqref{Ndelaskewcoop} is said to be {\em
uniformly stable for the order $\leq_A$ in $B_k$\/}  if for
every $\varepsilon>0$ there is a $\delta>0$, called the {\em
modulus of uniform stability\/}, such that, if $s\geq 0$ and
$\textsf{d}(u(s,\w_0,x_0),x)\leq \delta$ for certain $x\in B_k$
with $x\leq_A u(s,\w_0,x_0)$ or  $u(s,\w_0,x_0)\leq_A x$, then
for each $t\geq 0$,
\[ \di(u(t+s,\w_0,x_0),u(t,\w_0{\cdot}s,x))=\di(u(t,\w_0{\cdot}s,u(s,\w_0,x_0)),u(t,\w_0{\cdot}s,x))\leq
\varepsilon \,.\]
If this happens for each $k>0$, the forward orbit is said to be
{\em uniformly stable for the order $\leq_A$ in bounded sets\/}.}
\end{defi}

\begin{defi}
{\rm Given $k>0$, a positively $\tau$-invariant set $K\subset\Om\times
BU$ is {\em uniformly stable for the order $\leq_A$ in $B_k$\/}
if for any $\varepsilon>0$ there exists a $\delta>0$, called the
{\em modulus of uniform stability\/}, such that, if $(\w,x)\in K$,
$(\w,y)\in \Om\times B_k$ are such that $\di(x,y)< \delta$ with
$x\leq_A y$ or $y\leq_A x$, then $\di(u(t,\w,x),u(t,\w,y))\leq
\varepsilon$ for all $t\ge 0$. If this happens for each $k>0$, $K$
is said to be {\em uniformly stable for the order $\leq_A$ in
bounded sets\/}}.
\end{defi}

From Proposition~4.1 of~\cite{NOS2007} and Proposition~4.2
of~\cite{MNO2008}, if $z(t,\w_0,x_0)$ is bounded we  deduce that the family
$\{u(t,\w_0,x_0)\mid t\geq 0\}$ is equicontinuous on $[-k,0]$, and we can define
the omega-limit set of the trajectory of the point $(\w_0,x_0)$ as
\[\mathcal{O}(\w_0,x_0)=\{(\w,x)\in\Om\times BU\mid \exists
\,t_n\uparrow \infty\;\text{ with }
 \w_0{\cdot}t_n\to\w\,,\;
u(t_n,\w_0,x_0)\stackrel{\textsf{d}\;}\to x\}\,.\] Notice that the
omega-limit set of a pair $(\w_0,x_0)\in\Om\times BU$ makes sense
whenever $\cls_{X}\{u(t,\w_0,x_0)\mid t\geq 0\}$ is a
compact set.  As shown in Proposition~4.4
of~\cite{NOS2007} and~\cite{MNO2008}, $\mathcal{O}(\w_0,x_0)$ is a
positively invariant compact subset admitting a flow
extension.

Next statement shows that the omega-limit set inherits and improves the
stability properties of certain  relatively compact trajectories.

\begin{prop}\label{inherit}
Let $(\w_0,x_0)\in\Om\times BU$ with forward orbit
$\{\tau(t,\w_0,x_0)\,|\;t\geq 0\}$ relatively compact for the
product metric topology, and uniformly stable for $\leq_A$ in
bounded sets. Let $K$ denote the omega-limit set of $(\w_0,x_0)$.
Then

\begin{itemize}
\renewcommand{\labelenumi}{(\roman{enumi})}
\item[{\rm(i)}]   $K$ is uniformly stable for $\leq_A$ in bounded
sets.
\item[{\rm(ii)}]  $(K,\tau,\R^+)$ is uniformly stable.
\end{itemize}
\end{prop}

\begin{proof}
(i) Let $k_0>0$ be such that $\cls_{\Om\times
X}\{\tau(t,\w_0,x_0)\,|\;t\geq 0\}\subset \Om\times B_{k_0}\,$. Given
$k>0$, we check that $K$ is uniformly stable for $\leq_A$ in
$B_k$. Thus, we fix an $\ep>0$ and take $\delta>0$ the modulus of
uniform stability of the trajectory for $\leq_A$ in $B_{2k_0+k}$.
Let $(\w,x)\in K$ and $(\w,y)\in\Om\times B_k$ with
$\di(x,y)<\delta$ and $x\leq_A y$ or $y\leq_A x$, and take a
sequence $\{t_n\}_{n\in\N}$ such that
$\lim_{n\to\infty}\tau(t_n,\w_0,x_0)=(\w,x)$. Since
$u(t_n,\w_0,x_0)\leq_A u(t_n,\w_0,x_0)+y-x$ or
$u(t_n,\w_0,x_0)+y-x\leq_A u(t_n,\w_0,x_0)$,
$u(t_n,\w_0,x_0)+y-x\in B_{2k_0+k}$, and $\di(u(t_n,\w_0,x_0),
u(t_n,\w_0,x_0)+y-x)<\delta$ we deduce that for each $n\in\N$ and
$t\geq 0$
\[\di(u(t,\w_0{\cdot}t_n,u(t_n,\w_0,x_0)),u(t,\w_0{\cdot}t_n,u(t_n,\w_0,x_0)+y-x)))\leq
\ep\,.\] Finally, from Proposition 4.2 of~\cite{NOS2007}, and the
analogous one for the neutral case, which ensures certain continuity of the semiflow when the compact-open topology is considered on $BU$, as $n\uparrow \infty$ we get
$\di(u(t,\w,x),u(t,\w,y))\leq \varepsilon$ for each $t\geq 0$, and
(i) is proved.

(ii) Case I, that is, FDEs
and $\leq_A$ denotes  $\leq_{A,r}$. Let $(\w,x)$ and
$(\w,y)\in K$. Since $K$  is a positively invariant compact subset
admitting a flow extension, from Proposition~\ref{clasec1} we
deduce that $x,\,y\in C^1((-\infty,0],\R^m)$. We define
\begin{equation}\label{defiaxyr}
\begin{array}{ccl}
h\colon [-r,0] & \!\!\longrightarrow &\R^m\\
\quad s &\!\!\mapsto& \inf\{x'(s)-A\,x(s), y'(s)-A\,y(s)\}\,,
\text{
and} \\[.2cm]
a_{x,y}\colon  (-\infty,0]& \!\!\longrightarrow &\R^m\\\quad s
&\mapsto&
\begin{cases} i(s)=\inf\{x(s),y(s)\}\,, &
s\leq -r\\[.1cm]
e^{A\,(s+r)}i(-r)+\int_{-r}^s e^{A(s-\tau)} h(\tau)\,d\tau\,, &-r\leq s
\end{cases}
\end{array}\end{equation}
It is easy to check that $a_{x,y}\in BU$ and $a_{x,y}\leq x$. Moreover if
$-r\leq s\leq 0$,
\[ a'_{x,y}(s)=A\,a_{x,y}(s) +h(s)\,,\]
from which we deduce that $x'(s)- a'_{x,y}(s)\geq A(x(s)-a_{x,y}(s))$
and hence,  as stated in Remark~\ref{lipschitz}, $a_{x,y}\leq_A x$. Analogously,  $a_{x,y}\leq_A y$.

Moreover, since there is a $k_0>0$ such that
$\cls_{X} \{u(t,\w_0,x_0)\,|\;t\geq 0\}\subset B_{k_0}\,$,
$x'(s)=F(\w{\cdot}s,x_s)$ for each $s\in(-\infty,0]$ and $(\w,x)\in K$,
and there is $k>0$ such that $\n{F(\w,x)} \leq k  $ for each
$(\w,x)\in \Om\times B_{k_0}$, we deduce that there is a $k_1>0$ such
that $a_{x,y}\in B_{k_1}$ for each $(\w,x)$, $(\w,y)\in K$.

Let $\delta_1>0$ be the modulus of uniform
stability of $K$ for $\leq_A$ in $B_{k_1}$ for $\varepsilon/2$.
From the definition of $a_{x,y}$ and (F3), given $\delta_1>0$ we
can find $\delta>0$ such that
 if  $\di(x,y)\leq \delta$ then $\di(x,a_{x,y})\leq \delta_1$ and
$\di(y,a_{x,y})\leq \delta_1$. We omit the details, which are
similar to those of Cases II-III sketched below. Consequently,
whenever $\di(x,y)\leq \delta$, since $a_{x,y}\leq_A x$,
$a_{x,y}\leq_A y$ and $a_{x,y}\in B_{k_1}$ the uniform stability
of $K$ for the order $\leq_A$ yields
\[\di(u(t,\w,x),u(t,\w,y))\leq
\di(u(t,\w,x),u(t,\w,a_{x,y}))+\di(u(t,\w,a_{x,y}),u(t,\w,y))\leq
\ep\] and $(K,\tau,\R^+)$ is uniformly stable.

Cases II-III, i.e. FDEs and NFDEs when $\leq_A$ denotes the
ordering $\leq_{A,\infty}$. Both cases can be studied together by
taking $D\colon BU \to \R^m$, $x\mapsto x(0)$, in Case II. As in the
previous case, given $(\w,x)$ and $(\w,y)\in K$ we know that
$x,\,y\in C^1((-\infty,0],\R^m)$ and  from (A1) we can define
\begin{equation}\label{defiaxy}
\begin{array}{ccl}
a_{x,y}\colon(-\infty,0]&\longrightarrow &\R^m\\
\quad\;\; s &\mapsto& {\displaystyle \int_{-\infty}^s
e^{A(s-\tau)}\inf\{x'(\tau)-A\,x(\tau), y'(\tau)-A \,y(\tau)\}\,
d\tau \,,}
\end{array}\end{equation}
which satisfies $a_{x,y}\in BU$, $a_{x,y}\leq_A x$ and
$a_{x,y}\leq_A y$. Let $\widehat x=\widehat D x$ and $\widehat
y=\widehat D y$. As in Proposition~\ref{clasec1}
\begin{equation}\label{xyprima}
x'(s)=\int_{-\infty}^0 [d\widehat
\mu(\theta)]\,\widehat x\,'(s+\theta)\,,\quad y'(s)=\int_{-\infty}^0 [d\widehat \mu(\theta)]\,\widehat
y\,'(s+\theta)\,,
\end{equation}
and $\widehat x\,'(s)=G(\w{\cdot}s,\widehat x_s)$, $\widehat
y\,'(s)=G(\w{\cdot}s,\widehat y_s)$ for each
$s\in(-\infty,0]$.

Hence, since there is a $k_0>0$
such that $\cls_{X} \{u(t,\w_0,x_0)\,|\;t\geq 0\}\subset B_{k_0}\,$
and $k>0$ such that $\n{G(\w,\widehat x)} \leq k $ for each
$(\w,\widehat x)\in \Om\times B_{k_0\,\|\widehat D\|}$, we deduce
that there is a $k_1>0$ such that $a_{x,y}\in B_{k_1}$ for each
$(\w,x)$, $(\w,y)\in K$.

As in Case I, in order to
finish the proof it is enough to check that given $\delta_1>0$
there is a $\delta>0$ such that $\di(a_{x,y},x)\leq\delta_1$ and
$\di(a_{x,y},y)\leq\delta_1$ provided that $\di(x,y)\leq\delta$.
Assume on the contrary that there are $\delta_1>0$ and sequences
$\{(\w,x_n)\}_{n\in\N}$, $\{(\w,y_n)\}_{n\in\N}\subset K$ such
that $\lim_{\,n\to\infty}\di(x_n,y_n)=0$ and
$\di(a_{x_n,y_n},x_n)>\delta_1$. However, $\n{a_{x_n,y_n}(s)-x_n(s)}\leq \n{a_{x_n,y_n}(s)-b_{x_n,y_n}(s)}$ with
\[b_{x_n,y_n}(s)=\int_{-\infty}^s
e^{A(s-\tau)}\sup\{x_n'(\tau)-A\,x_n(\tau), y_n'(\tau)-A \,y_n(\tau)\}\,
d\tau \,,\]
from which we deduce that
\begin{equation}\label{axnyn-xn}
\n{a_{x_n,y_n}(s)-x_n(s)}\leq \int_{-\infty}^s\big\|e^{A (s-\tau)}\big\|
(\n{x_n'(\tau)-y_n'(\tau)}+\n{A}\,\n{x_n(\tau)-y_n(\tau)}) \,d\tau\,.
\end{equation}
Moreover, from relation~\eqref{xyprima} and
$G(\w,x)=F(\w,\widehat{D}^{-1}x)$ we have
\[\n{x_n'(\tau)-y'_n(\tau)}\leq \left \|\int_{-\infty}^0
[d \widehat
\mu(\theta)]\,(F(\w{\cdot}(\tau+\theta),(x_n)_{\tau+\theta})
-F(\w{\cdot}(\tau+\theta),(y_n)_{\tau+\theta}))\,\right\|\,.
\]
Hence, from the total finite variation of $\widehat \mu$,
hypothesis (F2), the uniform continuity of $F$ on $K$, implied
by hypothesis (F3), hypothesis (A1) and $\lim_{\,n\to\infty}\di(x_n,y_n)=0$, we
deduce that $a_{x_n,y_n}-x_n$ tends to $0$ uniformly on compact
sets of $(-\infty,0]$,
 which contradicts that $\di(a_{x_n,y_n},x_n)>\delta_1$ and finishes the
proof for Cases II-III.
\qquad\end{proof}

The next result shows that, under the previous assumptions, the
omega-limit set $K$ is always a  minimal subset in Case I,
and more regularity in the initial data is needed for Cases II and
III.

First we recall the definition of the section map of a compact set
$M\subset \Om\times X$. We introduce the projection set of $M$
into the fiber~space
\[
M_X=\{x\in X \mid \text{ there exists } \w\in \Om \text{ such that
} (\w,x)\in M \}\subseteq X\,.
\]
From the compactness of $M$ it is immediate to show that also
$M_X$ is a compact subset of $X$. Let $\mathcal P_c(M_{X})$ denote
the set of closed subsets of $M_X$, endowed with the Hausdorff
metric $\rho$, that is, for any two sets $C$, $B\in \mathcal
P_c(M_X)$,
\[
\rho(C,B)=\sup\{\alpha(C,B),\;\alpha(B,C)\}\,,
\]
where $\alpha(C,B)=\sup\{r(c,B)\mid c\in C\}$ and
$r(c,B)=\inf\{\di(c,b)\mid b\in B\}$. Then, define the so-called
{\em section map\/}
\begin{equation}\label{sectionmap}
\begin{array}{ccl}
  \Om  & \longrightarrow & \mathcal P_c(M_X)\\
    \w & \mapsto & M_\w=\{x\in
X \mid (\w,x)\in M \}\,.
\end{array}
\end{equation}
Due to the minimality of $\Om$ and the compactness of $M$, the set
$M_\w$ is nonempty for every $\w\in \Om$; besides, the map is
trivially well-defined.

\begin{teor}\label{minimal}
Let $(\w_0,x_0)\in\Om\times BU$ with forward orbit
$\{\tau(t,\w_0,x_0)\,|\;t\geq 0\}$ relatively compact for the
product metric topology, and uniformly stable for $\leq_A$ in
bounded sets. Let $K$ denote the omega-limit set of $(\w_0,x_0)$.
Then, $K$ is a minimal subset in {\rm Case I}, and
in {\rm Cases II} and {\rm III} provided that $x_0$ is
Lipschitz.
\end{teor}

\begin{proof}
Let $M$ be a minimal subset such that $M\subseteq K$.  We just
need to show that $K \subseteq M$. So, take an element $(\w,x)\in
K$ and let us prove that $(\w,x)\in M$. As $M$ is in particular
closed, it suffices to see that for any fixed $\ep>0$ there exists
$(\w,x^*)\in M$ such that $\di(x,x^*)\leq \ep$.

First of all, there exists $s_n\uparrow\infty$ such that
$\lim_{\,n\to\infty}(\w_0{\cdot}s_n,u(s_n,\w_0, x_0))=(\w,x)$. Now, take
a pair $(\w, \wt x)\in M\subseteq K$. Then, there exists a
sequence $t_n\uparrow \infty$ such that
\[
(\w, \wt x)=\lim_{n\to \infty} (\w_0{\cdot}t_n,u(t_n,\w_0,x_0))\,.
\]
Since from Proposition~\ref{inherit} $(M,\tau,\R^+)$  is uniformly
stable, we can apply Theorem~3.4 of~\cite{NOS2007} so that the
section map~\eqref{sectionmap} turns out to be continuous at any
point. As $\w_0{\cdot}t_n\to \w$, we deduce  that $M_{\w_0\cdot t_n}\to
M_{\w}$ in the Hausdorff metric. Therefore, for $\wt x\in M_{\w}$
there exists a sequence $x_n\in M_{\w_0\cdot t_n}$, $n\geq 1$,
such that $x_n\to \wt x$ as $n\uparrow \infty$. From
Proposition~\ref{clasec1} we deduce that $x_n\in C^1((-\infty,0],\R^m)$
for each $n\in\N$ and denoting $y_n=u(t_n,\w_0,x_0)$, we have
$\di(x_n,y_n)\to 0$ as $n\uparrow\infty$. The rest of the proof
depends on the case we are dealing with.

Case I, that is, FDE and $\leq_A$ denotes
$\leq_{A,r}$. Let $n_0$ be such that $t_n> r$ for each $n\geq
n_0$. Then $y_n\in C^1([-r,0],\R^m)$ and we can define $a_{x_n,y_n}$ as
in~\eqref{defiaxyr} for each $n\geq n_0$. Moreover, as in
Proposition~\ref{inherit}, we check that $a_{x_n,y_n}\in BU$,
$a_{x_n,y_n}\leq_A x_n$, $a_{x_n,y_n}\leq_A y_n$ and there is
a  $k_1>0$ such that $a_{x_n,y_n}\in B_{k_1}$ for each $n\geq
n_0$.

Next, let $\delta>0$ be the modulus of uniform stability of $K$
for $\leq_A$ in $B_{k_1}$ for $\varepsilon/2$. Since
$\lim_{\,n\to\infty}\di(x_n,y_n)=0$, from the definition of
$a_{x_n,y_n}$ and hypothesis (F3) there is an $n_1\geq n_0$ such that
$\di(x_n,a_{x_n,y_n})<\delta$ and $\di(y_n,a_{x_n,y_n})<\delta$
for each $n\geq n_1$. Hence, the uniform stability of $K$ for the
order $\leq_A$ in $B_{k_1}$ yields
\[\di(u(t+t_{n_1},\w_0,x_0),u(t,\w_0{\cdot}t_{n_1},x_{n_1}))=
\di(u(t,\w_0{\cdot}t_{n_1},y_{n_1}),u(t,\w_0{\cdot}t_{n_1},x_{n_1}))\leq \ep\] for each $t\geq 0$. In particular, if $n_2$ is such that
$s_n-t_{n_1}\geq 0$ for $n\geq n_2$, we obtain
\begin{equation}\label{uno}
\di(u(s_n,\w_0,x_0),u(s_n-t_{n_1},\w_0{\cdot}t_{n_1},x_{n_1}))\leq \ep
\quad \text{for each } \; n\geq n_2\,.
\end{equation}
Now, it remains to notice that, as $(\w_0{\cdot}t_{n_1},x_{n_1})\in M$,
also
$\tau(s_n-t_{n_1},\w_0{\cdot}t_{n_1},x_{n_1})=(\w_0{\cdot}s_n,u(s_n-t_{n_1},\w_0{\cdot}t_{n_1},x_{n_1}))\in
M$ for all $n\geq n_2$. Therefore, there is a convergent
subsequence towards a pair $(\w,x^*)\in M$, and taking limits
in~\eqref{uno}, we deduce that $\di(x,x^*)\leq \ep$, as we wanted.

Cases II-III, that is, FDEs and NFDEs when $\leq_A$ denotes the
ordering $\leq_{A,\infty}$. Both cases can be studied together by
taking $D\colon BU \to \R^m$, $x\mapsto x(0)$, in Case II. Remember
that $x_n\in C^1((-\infty,0],\R^m)$, and since $x_0$ is Lipschitz,
it is not hard to check that $y_n=u(t_n,\w_0,x_0)$ is also
Lipschitz. Therefore,
\[\inf\{x_n'(\tau)-A \,x_n(\tau), y_n'(\tau)-A \,y_n(\tau)\}\]
is defined for almost every
$\tau\in(-\infty,0]$ and we can define $a_{x_n,y_n}$ as
in~\eqref{defiaxy}. As in Proposition~\ref{inherit}, we check that
$a_{x_n,y_n}\in BU$, $a_{x_n,y_n}\leq_A x_n$,
$a_{x_n,y_n}\leq_A y_n$ and there is a  $k_1>0$ such that
$a_{x_n,y_n}\in B_{k_1}$ for each $n\in \N$.

Moreover from relation~\eqref{usombrero},  $y_n= \widehat{D}^{-1}
\widehat u(t_n,\w_0, \widehat D x_0)=\widehat{D}^{-1} \widehat
u(0,\w_0{\cdot}t_n, \widehat D y_n)$ and
\[(\widehat Dy_n)'(s)=G(\w_0{\cdot}(t_n+s), (\widehat Dy_n)_s)\,,\quad -t_n\leq s\leq 0\,.\]
Therefore, we deduce  from~\eqref{defiT} that
\begin{equation*}
 y_n(s)=\int_{-\infty}^0
[d\widehat \mu(\theta)]\,\widehat D y_n(s+\theta)\,,\quad s\in
(-\infty,0]\,,
\end{equation*}
which yields
\begin{multline}\label{increyn}
\frac{y_n(s+h)-y_n(s)}{h}= \int_{-\infty}^{-t_n-s} [d\widehat \mu(\theta)]\left[ \frac{\widehat Dy_n(s+\theta+h)-\widehat Dy_n(s+\theta) }{h}\right]\\ +
 \int_{-t_n-s}^0
[d\widehat \mu(\theta)]\left[\frac{1}{h}\int_{s+\theta}^{s+\theta+h} G(\w_0{\cdot}(t_n+t),(\widehat D
y_n)_t)\,dt \right]\,,
\end{multline}
for $-t_n< s\leq s+h\leq 0$.
Analogously, from
\[\frac{d}{ds}D(x_n)_s=F(\w_0{\cdot}(t_n+s), (x_n)_s)\,,\quad s\in(-\infty,0]\,\]
and
$G(\w,x)=F(\w,\widehat{D}^{-1}x)$ we deduce that
\begin{multline}\label{increxn}
\frac{x_n(s+h)-x_n(s)}{h}= \int_{-\infty}^{-t_n-s} [d\widehat \mu(\theta)]\left[ \frac{\widehat D x_n(s+\theta+h)- \widehat Dx_n(s+\theta)}{h}\right]\\ + \int_{-t_n-s}^0 [d\widehat \mu(\theta)]
\left[\frac{1}{h}\int_{s+\theta}^{s+\theta+h} G(\w_0{\cdot}(t_n+t),(\widehat D
x_n)_t)\,dt \right] \,,
\end{multline}
for $-t_n< s\leq s+h\leq 0$. Moreover, since $x'_n$ and $y'_n$ exist almost everywhere and $\widehat D y_n$ and $\widehat D x_n$ are Lipschitz on $(-\infty,0]$ with the same constant for all $n\in\N$, from
relations \eqref{increyn}, \eqref{increxn} and
$G(\w,x)=F(\w,\widehat{D}^{-1}x)$ we deduce that there is a constant $c_0\geq 0$ such that
\begin{multline*}
\n{x_n'(\tau)-y'_n(\tau)}\leq c_0\,\n{\widehat \mu}_\infty(-\infty,-t_n-\tau]\\+ \left\| \int_{-t_n-\tau}^0 \!\![d\widehat \mu (\theta)]\,(F(\w_0{\cdot}(t_n+\tau+\theta),(x_n)_{\tau+\theta}) - F(\w_0{\cdot}(t_n+\tau+\theta),(y_n)_{\tau+\theta}))\right\|
\end{multline*}
 for almost every
$\tau\geq -t_n$. Hence, from~\eqref{axnyn-xn},  the total
finite variation of $\widehat \mu$, hypotheses (F2), the
uniform continuity of $F$ in $\cls_{\Om\times
X}\{\tau(t,\w_0,x_0)\,|\;t\geq 0\}$, implied by hypothesis (F3),
$t_n\uparrow \infty$, hypothesis (A1), and $\lim_{\,n\to\infty}\di(x_n,y_n)=0$ we
deduce that
\[\lim_{\,n\to\infty} \di(x_n,a_{x_n,y_n})=0\quad \text{ and
}\quad \lim_{\,n\to\infty} \di(y_n,a_{x_n,y_n})=0\,.
\]
Consequently,  if $\delta>0$ is the modulus of uniform stability
of $K$ for $\leq_A$ in $B_{k_1}$ for $\varepsilon/2$, there is
an $n_1$ such that $\di(x_n,a_{x_n,y_n})<\delta$ and
$\di(y_n,a_{x_n,y_n})<\delta$ for each $n\geq n_1$, and the rest
of the proof is identical to that of Case I.
\qquad\end{proof}

\section{Topological structure of omega-limit sets}\label{copiadelabase}
As in the previous section we consider the
monotone skew-product semiflows~\eqref{Ndelaskewcoop} induced by~\ref{Ninfdelay} in Cases
I-III. First,  we extend to this setting results
of Novo, N\'{u}\~{n}ez and Obaya~\cite{NNO2005} and Novo, Obaya and Sanz~\cite{NOS2007} ensuring the presence of almost
automorphic dynamics from the existence of a semicontinuous
semi-equilibrium.

\begin{defi}\label{equilibrios}
{\rm A map $a:\Om\to BU$ such that $u(t,\w,a(\w))$ is defined for any
$\w\in \Om$, $t\geq 0$~is

\begin{itemize}
\item[a)] an {\em equilibrium\/} if $a(\w{\cdot}t)=u(t,\w,a(\w))$ for
any $\w\in\Om$ and $t\ge 0$, \item[b)] a {\em super-equilibrium\/}
if $a(\w{\cdot}t)\ge_A u(t,\w,a(\w))$ for any $\w\in\Om$ and $t\ge 0$,
and
\item[c)] a {\em sub-equilibrium\/} if $a(\w{\cdot}t)\le_A u(t,\w,a(\w))$
for any $\w\in\Om$ and $t\ge 0$.
\end{itemize}
We will call {\em semi-equilibrium} to either a super or a
sub-equilibrium.}
\end{defi}

\begin{defi}\label{semicont}
{\rm A super-equilibrium (resp.~sub-equilibrium) $a:\Om\to BU$ is {\em
semicontinuous\/} if the following properties hold:

\begin{itemize}
\item[(1)] $\G_a=\cls_{X}\{a(\w)\,|\;\w\in\Om\}$ is a
compact subset of $X$ for the compact-open topology, and
\item[(2)] $C_a=\{(\w,x)\mid x\le_A a(\w)\}$
(resp.~$C_a=\{(\w,x)\mid x\ge_A a(\w)\}$) is a closed subset of
$\Om\times X$ for the product metric topology.
\end{itemize}
An equilibrium is {\em semicontinuous} in any of these cases.}
\end{defi}

As shown in Proposition 4.8 of~\cite{NOS2007} a semicontinuous
equilibrium does always have a residual subset of continuity
points. This theory requires topological properties of semicontinuous maps
stated in Aubin and Frankowska~\cite{aufr} and Choquet~\cite{choq}. The next result shows that a semicontinuous semi-equilibrium  provides an almost automorphic
extension of the base if a relatively compact trajectory exists. We omit its proof, analogous to the one
of Proposition 4.9 of~\cite{NOS2007} once Proposition~\ref{inherit} is proved.

\begin{prop}\label{almostauto}
Let $a:\Om\to BU$ be a semicontinuous semi-equilibrium and assume
that there is an $\w_0\in \Om$ such that  $\cls_{X}\{u(t,\w_0,a(\w_0))\mid t\ge 0\}$ is a compact subset
of $X$ for the compact-open topology. Then:
\begin{itemize}
\item[{\rm(i)}] The omega-limit set $\mathcal{O}(\w_0,a(\w_0))$
contains a unique minimal set, which is an almost automorphic
extension of the base flow. \item[{\rm (ii)}] If  the orbit
$\{\tau(t,\w_0,a(\w_0))\mid t\geq 0\}$ is uniformly stable for
$\leq_A$ in bounded sets, then $\mathcal{O}(\w_0,a(\w_0))$ is a
copy of the base.
\end{itemize}
\end{prop}

\vspace{.1cm}
If the semicontinuous semi-equilibrium satisfies some supplementary
and somehow natural compactness conditions, a semicontinuous
equilibrium is obtained. As shown in Proposition 4.10
of~\cite{NOS2007} for infinite delay, provided that $\G_a\subset BU$ and $\sup_{\w\in\Om} \n{a(\w)}_\infty<\infty$, it can be proved that the
following conditions are equivalent:

\vspace{.1cm}
\begin{itemize}
\item[(C1)]  $\G=\cls_{X}\{u(t,\w,a(\w))\,\mid t\ge
0,\,\w\in\Om\}$ is a compact subset of $BU$ for the compact-open
topology.
\item[(C2)] For each $\w\in\Om$, the $\cls_{X}
\{u(t,\w,a(\w)) \mid t\ge 0\}$ is a compact subset of $BU$ for the
compact-open topology.
\item[(C3)] There is  an $\w_0\in\Om$ such that
the $\cls_{X}\{u(t,\w_0,a(\w_0))\mid t\ge 0\}$ is a
compact subset of $BU$ for the compact-open topology.
\end{itemize}
Consequently, an easy adaptation of the proof of Theorem 4.11
of~\cite{NOS2007} proves the following result.

\begin{teor}\label{main}
Let us assume the existence of a semicontinuous semi-equilibrium
$a:\Om\to BU$ satisfying   $\sup_{\w\in\Om}
\n{a(\w)}_\infty<\infty$, $\G_a\subset BU$  and one of the
equivalent conditions {\rm(C1)-(C3)}. Then,

\begin{itemize}
\item[\rm(i)] there exists a semicontinuous equilibrium
$c:\Om\to BU$ with $c(\w)\in\G$ for any $\w\in\Om$. \item[\rm
(ii)] Let $\w_1$ be a continuity point for $c$. Then, the
restriction of the semiflow $\tau$ to the minimal set
\begin{equation*}
 K^*=\cls_{\Om\times X}\{(\w_1{\cdot}t,c(\w_1{\cdot}t))\,|\;t\ge
 0\}\subset C_a
\end{equation*}
is an almost automorphic extension of the base flow
$(\Om,\sigma,\R)$. \item[\rm(iii)] $K^*$ is the only minimal set
contained in the omega-limit set
$\mathcal{O}(\widehat\w,a(\widehat\w))$ for each point
$\widehat\w\in\Om$. \item[\rm(iv)] If there is a point
$\wt\w\in\Om$ such that the trajectory
$\{\tau(t,\wt\w,a(\wt\w))\mid t\geq 0\}$ is uniformly stable for
$\leq_A$ in bounded sets, then for each $\widehat\w\in\Om$,
\[\mathcal{O}(\widehat\w,a(\widehat\w))=K^*=\{(\w,c(\w))\mid\w\in\Om\}\,,\]
i.e.~it is a copy of the base determined by the equilibrium $c$ of
{\rm(i)}, which is a continuous map.
\end{itemize}
\end{teor}

For the rest of the section we will assume the following stability
assumption for the trajectories with initial data in a ball.

\vspace{.1cm}
\begin{itemize}
\item[(F5)] There is a $k_0>0$ such that all the trajectories with Lipschitz
initial data in $B_{\widehat k_0}$  are relatively compact  for the product metric
topology and uniformly stable for $\leq_A$ in bounded sets, where
\begin{equation}\label{k0sombrero}
\widehat k_0=\frac{\|{\widehat D}^{-1}\|\,
\sup\{\n{F(\w,x)}\mid (\w,x)\in \Om\times {\widehat D}B_{k_0}\}}{\n{A}}+k_0\,.
\end{equation}
\end{itemize}
The following result provides a continuous super-equilibrium for every compact, positively invariant set contained in $\Om\times B_{k_0}$. It requires basic properties of the exponential ordering obtained in Section~\ref{omegalimit}.

\begin{teor}\label{inferior}
Let $(\w_0,x_0)\in \Om\times BU$ with forward orbit
$\{\tau(t,\w_0,x_0)\,|\;t\geq 0\}$ relatively compact for the
product metric topology, uniformly stable for $\leq_A$ in
bounded sets and $\cls_{X} \{u(t,\w_0,x_0)\,|\;t\geq 0\}\subset B_{k_0}\,$. Let $K=\mathcal{O}(\w_0,x_0)\subset \Om\times
B_{k_0}$ be its omega-limit set. For each $\w\in\Om$ we define the
map $a(\w)$ on $(-\infty,0]$ by
\begin{align*}
&\hspace{-1cm} \text{{\rm Case I}}: \quad a(\w)(s)=\begin{cases}
i(\w)(s)=\inf\{x(s) \mid (\w,x)\in K\}\,, &
s\leq -r\,,\\[.1cm]
e^{A\,(s+r)}i(\w)(-r)+  \int_{-r}^s e^{A(s-\tau)}\,h(\w)(\tau)\,d\tau\,, & -r\leq s\,,
\end{cases}\\[.1cm]
&\hspace{-1cm}\text{{\rm Cases II-III}}: \quad a(\w)(s)= \int_{-\infty}^s e^{A(s-\tau)}\,h(\w)(\tau)\,d\tau\,,\quad
s\leq 0\,, \\[.2cm]
\text{ where }  & \quad \begin{array}[t]{ccl}
h(\w)\colon (-\infty,0] & \!\!\longrightarrow &\R^m\\
\quad \tau &\!\!\mapsto& \inf\{x'(\tau)-A\, x(\tau) \mid
(\w,x)\in K\}\,.\end{array}
\end{align*}
Then, $a(\w)$ is Lipschitz for every $\w\in\Om$,  $\sup_{\w\in\Om}
\n{a(\w)}_\infty\leq  \widehat k_0$, and the map $a\colon\Om\rightarrow BU$, $\w\mapsto a(\w)$ is well-defined, it is a continuous super-equilibrium and it satisfies the equivalent
conditions {\rm(C1)-(C3)}.
\end{teor}

\begin{proof}
Since $K$  is a positively invariant compact subset admitting a
flow extension, from Proposition~\ref{clasec1} we deduce that
$x\in C^1((-\infty,0],\R^m)$ for each $(\w,x)\in K$.
We claim that there is a positive constant $L>0$ such that $x$ is
a Lipschitz function with constant $L$ for each
$(\w,x)\in K$. In Case I and Case II this assertion follows~from
\[z'(t,\w_0,x_0)=F(\w_0{\cdot}t,u(t,\w_0,x_0))\,,\quad t\geq 0\,,\]
and (F2) with $L=\sup\{\n{F(\w,x)}\mid (\w,x)\in \Om\times B_{k_0}\}$.
In Case III, from relation~\eqref{usombrero}
\[u(t,\w_0,x_0)= {\widehat D}^{-1} {\widehat u}(t, \w_0, {\widehat D}x_0)\,,\]
and we deduce the result with $L=\|{\widehat D}^{-1}\|\,
\sup\{\n{F(\w,x)}\mid (\w,x)\in \Om\times {\widehat D}B_{k_0}\}$. From
this, as in~\cite{NOS2007} we deduce that $i(\w)$ is also
Lipschitz on $(-\infty,-r]$ with the same constant.

Since $\n{x}_\infty\leq k_0$ and $\n{x'}_\infty\leq
L$ for each $x\in K_\w$, the function $h(\w)$ is well defined for
each $\w\in\Om$. Moreover, since the family $\{x'(\tau)-A\,x(\tau)\mid (\w,x)\in K\}$
is equicontinous on $[-k,0]$, we can check that $h(\w)\in
X=C((-\infty,0],\R^m)$ and $\sup_{\w\in\Om}
\n{h(\w)}_\infty\leq L+\n{A}\,k_0$.

Then, $a(\w)$ is well
defined for each $\w\in\Om$, $a(\w)\in X$ and $\sup_{\w\in\Om}
\n{a(\w)}_\infty< \widehat k_0$, for
$\widehat k_0=k_0+L/\n{A}$, as defined in~\eqref{k0sombrero}.
 In addition, $a(\w)$ belongs to
$C^1([-r,0],\R^m)$ in Case I, and to $C^1((-\infty,0],\R^m)$ in Cases
II-III with
\[a(\w)'(s)=A\, a(\w)(s)+h(\w)(s)\,, \qquad\begin{array}{rl}
-r\leq s\leq 0\,, & \text{Case I}\,.\\
-\infty<s\leq 0\,, & \text{Cases II-III}\,.
\end{array}\]
From this fact, the uniform boundedness of $a(\w)$ and $h(\w)$,
and the uniform Lipschitz character of $i(\w)$ on $(-\infty,-r]$,
we deduce that there is a positive constant $\widehat L>0$ such
that $a(\w)$ is a Lipschitz function with constant
$\widehat L$ for each $\w\in \Om$. Hence, $a(\w)$ belongs to $BU$,
i.e. $a$ is well defined. Moreover, $\G_a=\cls_{X}\{a(\w)\,|\;\w\in\Om\}$
 is a compact subset of $X$, and actually $\G_a\subset BU$.

Let us check that $a$ defines a super-equilibrium.  Since
$a(\w)\in B_{\widehat k_0}$, it follows from hypothesis (F5) that $u(t,\w,a(\w))$ exists
for any $\w\in \Om$ and $t\geq 0$. It is easy to prove, as in
Proposition~\ref{inherit}, that $a(\w)\leq_A x$  for each $x\in
K_\w\,$. Next we claim that if $z\in BU$ with $z\leq_A x$ for
each $x\in K_\w$ then $z\leq_A a(\w)$ provided that $z$ is
Lipschitz on $[-r,0]$ in Case I, and on $(-\infty,0]$ in
Cases~II-III.

From the definition of $a(\w)$ it is
easy to check that $z(s)\leq
a(\w)(s)$ for each $s\leq 0$. Moreover, since $x\in C^1((-\infty,0],\R^m)$, $z\leq_A x$ and $z$ is Lipschitz we deduce from Remark~\ref{lipschitz} that
\[z'(s)-A\, z(s)\leq x'(s)-A\, x(s)\]
for almost every $s\in[-r,0]$ in Case I, almost every
$s\in(-\infty,0]$ in Cases II-III, and every $x\in K_\w$. Hence,
the definition of $h(\w)$ provides at these points
\[z'(s)-A\, z(s)\leq h(\w)(s)=a(\w)'(s)-A\,a(\w)(s)\,,\]
and we conclude that $z\leq_A a(\w)$,
as claimed.

Fix $\w\in \Om$, $t\geq 0$ and consider
any $y\in K_{\w{\cdot}t}$, i.e. $(\w{\cdot}t,y)\in K$. As we have a flow on
$K$, $\tau(-t,\w{\cdot}t,y)=(\w,u(-t,\w{\cdot}t,y))\in K$ and therefore,
$a(\w)\leq_A u(-t,\w{\cdot}t,y)$. Applying monotonicity,
$u(t,\w,a(\w))\leq_A y$. As this happens for any $y\in
K_{\w{\cdot}t}$, and $u(t,\w,a(\w))$ is Lipschitz on $(-\infty,0]$ we
get that $u(t,\w,a(\w))\leq_A a(\w{\cdot}t)$ and $a$ is a
super-equilibrium, as stated.

Now let us prove that $a$ is continuous on $\Om$. From the
definition of $a$ it is enough to check the continuity of $h$ on
$\Om$, as well as the continuity of $i$ in Case~I. The continuity
of $i$ is shown in Proposition 5.2 of~\cite{NOS2007}, so let us
prove that $h$ is continuous on $\Om$. Fix $\w\in\Om$ and assume
that  $\w_n\to\w$ and $h(\w_n)\stackrel{\textsf{d}\;}\to y$ as
$n\uparrow \infty$. First we check that $h(\w)\leq y$. Fix $s\in
(-\infty,0]$ and $i\in\{1,\ldots,m\}$.
From the definition of $h$ there are $(\w_n,x_n)\in
K$, depending on $s$ and $i$, although dropped from the notation, such that
\[\left|h_i(\w_n)(s)-(x'_{n,i}(s) -(A \,x_n)_i(s))\right|<\frac{1}{n}\,,\]
 where  $x_{n,i}$ and $(A\,x_n)_i$  denote the corresponding components of $x_n$ and $A\,x_n$. This implies that
$\lim_{\,n\to\infty}(x'_{n,_i}(s)- (A\, x_n)_i(s))=y_i(s)$ for $i\in\{1,\ldots,m\}$, that is,
\[\lim_{n\to\infty} x'_n(s)- A\,x_n(s)=y(s)\,.\]
Moreover, from the
compactness of $K$, an adequate subsequence $(\w_{n_j},x_{n_j})$
tends to some $(\w,x)\in K$ in the product metric topology. From
this, it can be shown that $\lim_{\,n\to\infty}x'_{n_j}(s)=x'(s)$;
hence $y(s)=x'(s)-A\,x(s)$, and again, from the definition of
$h(\w)$ we conclude that $h(\w)(s)\leq y(s)$. As this happens for
each $s\in(-\infty,0]$, we deduce that $h(\w)\leq y$. On the
other hand, from Proposition~\ref{inherit} we know that
$(K,\tau,\R^+)$ is uniformly stable, and then Theorem~3.4
of~\cite{NOS2007} asserts that the section map for $K$, $\w\in\Om
\mapsto K_\w$, is continuous at every $\w\in\Om$, which   implies
that  $K_{\w_n}\to K_\w$ in the Hausdorff metric. Therefore, for
any $z\in K_\w$ there exist $z_n\in K_{\w_n}$, $n\geq 1$, such
that $z_n \stackrel{\textsf{d}\;}\to z$. Then, $(\w_n,z_n)\in K$
implies that $h(\w_n)(s)\leq z_n'(s)- A\, z_n(s)$ and taking
limits, $y(s)\leq z'(s)-A\, z(s)$ for each $s\in(-\infty,0])$.
As this happens for any $z\in K_\w$, we conclude that $y\leq
h(\w)$. In all, $h(\w)=y$, as wanted. Consequently
$\G_a=\{a(\w)\,|\;\w\in\Om\}$.

Finally,
since $a(\w)\in B_{\widehat k_0}$ for each $\w\in\Om$, from hypothesis
(F5), Proposition~4.1 of~\cite{NOS2007} and Proposition~4.2
of~\cite{MNO2008}, we deduce that the equivalent conditions
(C1)-(C3) hold, and the proof is complete. Notice that $a(\w)$ is the infimum
among the Lipschitz functions of the set $K_\w$.
\qquad\end{proof}

For the main theorem of the paper, in which the 1-covering property of omega-limit sets
is established, we add the following hypothesis:

\begin{itemize}
\item[(F6)]
   If $(\w,x)$, $(\w,y)\in\Om\times BU$ admit a backward orbit extension,  $x\leq_A y$, and there is a subset $J\subset\{1,\ldots,m\}$ such that
\begin{equation*}
\begin{split}
x_i=y_i & \quad \text{ for each } i\notin J\,,\\
x_i(s)< y_i(s) & \quad \text{ for each } i\in J\;\text{ and } s\leq 0\,
\end{split}
\end{equation*}
    then $F_i(\w,y)-F_i(\w,x) - (A\,(Dy-Dx))_i>
0$ for each $i\in J$ and $\w\in\Om$\,,
\end{itemize}
for all the cases, and the following ones for Case III:

\begin{itemize}
\item[(D6)] The measure $\widehat\mu$ of~\eqref{defiT} is positive, i.e. $\widehat D^{-1}$ is a positive operator.
\item[(D7)] $AD\,x=DA\,x$ for each $x\in BU$.
\end{itemize}

\begin{teor}\label{maincopy} We consider the
monotone skew-product semiflows~\eqref{Ndelaskewcoop} induced by~{\rm\ref{Ninfdelay}} in the cases

{\rm \vspace{.2cm}\begin{tabular}{ll}
        \!\!Case I: & FDEs~\ref{infdelay} under (F1)-(F6) and order  $\leq_{A,r}$;\\
        \!\!Case II: & FDEs~\ref{infdelay} under (F1)-(F6), (A1) and order  $\leq_{A,\infty}$;\\
        \!\!Case III: & NFDEs~\ref{Ninfdelay}, under (F1)-(F6), (A1), (D1)-(D7) and order  $\leq_{A,\infty}$;\vspace{.2cm}
\end{tabular}}

\noindent with $D$ of the form~\eqref{Dformula}, $\nu\equiv 0$ in {\rm Cases I-II}.
Let $(\w_0,x_0)\in \Om\times B_{k_0}$ with forward orbit
$\{\tau(t,\w_0,x_0)\,|\;t\geq 0\}$ relatively compact for the
product metric topology, and uniformly stable for $\leq_A$ in
bounded sets, be such that its omega-limit set $K=\mathcal{O}(\w_0,x_0)\subset
\Om\times B_{k_0}$. In addition, in {\rm Cases II} and {\rm III}
assume that $x_0$ is Lipschitz.

Then $K=\mathcal{O}(\w_0,x_0)=\{(\w,c(\w))\mid \w\in\Om\}$ is a copy of
the base and {\rm \[\lim_{t\to\infty}
\di(u(t,\w_0,x_0),c(\w_0{\cdot}t))=0\,,\]} where $c:\Om\to BU$ is a
continuous equilibrium.
\end{teor}

\begin{proof} We apply
Theorem~\ref{inferior} to obtain  a continuous super-equilibrium  $a$ satisfying (C1)-(C3)
with $a(\w)$ Lipschitz for each $\w\in\Om$. Then, from (F5) and Theorem~\ref{main} we deduce that
there is a continuous equilibrium $c:\Om\to BU$ such that for each~$\widehat\w\in\Om$,
\begin{equation}\label{copia}
\mathcal{O}(\widehat\w,a(\widehat\w))=K^*=\{(\w,c(\w))\mid
\w\in\Om\}\,.
\end{equation}
The definition of $a$ yields  $a(\w)\leq_A x$ for each
$(\w,x)\in K$ and hence  $c(\w)\leq_A x$ by the construction of
$c$. As in
Jiang and Zhao~\cite{jizh}  and Novo, Obaya and Sanz~\cite{NOS2007} we prove that there is a subset
$J\subset\{1,\ldots,m\}$ such that
\begin{equation*}
\begin{split}
c_i(\w)=x_i & \quad \text{ for each } (\w,x)\in K \text{ and  }i\notin J\,,\\
c_i(\w)< x_i & \quad \text{ for each } (\w,x)\in K \text{ and
}i\in J\,.
\end{split}
\end{equation*}
It is enough to check that if $c_i(\wt\w)(0)=\wt x_i(0)$ for some
$i\in\{1,\ldots,m\}$ and $(\wt\w,\wt x)\in K$, then $c_i(\w)=x_i$
for any $(\w,x)\in K$. We first notice that $c_i(\wt\w)=\wt x_i$.
Otherwise, there would be an $s_0\in(-\infty,0]$ with
$c_i(\wt\w)(s_0)<\wt x_i(s_0)$.

Cases II-III: from
$c(\wt\w)\leq_A \wt x$ we know that
\[\wt x(0)-c(\wt
\w)(0)\geq e^{-As_0}(x(s_0)-c(\wt\w)(s_0))\]
 which
implies from Remark~\ref{quasipositive} that $c_i(\wt\w)(0)<\wt x_i(0)$, a contradiction.

 Case I: in this case $\leq_A$ denotes $\leq_{A,r}$. Since  $K$ admits a flow extension, for each $s\leq 0$ we know that
 $(\wt \w{\cdot}s,u(s,\wt\w,\wt x))=(\wt\w{\cdot}s,\wt x_s)\in K$. Hence,
$c(\wt\w{\cdot}s)\leq_A \wt x_s$  and
\[\wt x_s(t)-c(\wt \w{\cdot}s)(t)\geq e^{A (t-\theta)}(\wt x_s(\theta)-c(\wt \w{\cdot}s)(\theta))\,,\quad -r\leq \theta\leq t\leq 0\,,\; s\leq 0\,.\]
Therefore, as above, from $c_i(\wt\w)(s_0)<\wt x_i(s_0)$ and Remark~\ref{quasipositive}, in a finite number of steps, we deduce that $c_i(\wt\w)(0)<\wt x_i(0)$,
a contradiction.

Therefore, $c_i(\wt \w)=\wt x_i$. Next,
from Theorem~\ref{minimal} we know that $K$ is minimal.
Thus we take $(\w,x)\in K$ and a sequence $s_n\downarrow -\infty$
such that $\wt\w{\cdot}s_n\to\w$ and $u(s_n,\wt \w,\wt
x)\stackrel{\textsf{d}\;}\to x$.~Then,
\begin{align*}
     x_i(0) & =\lim_{n\to\infty} u_i(s_n,\wt\w,\wt x)(0)=
        \lim_{n\to\infty} \wt x_i(s_n) \\  &
        =\lim_{n\to\infty} c_i(\wt\w)(s_n)
        =\lim_{n\to\infty}c_i(\wt\w{\cdot}s_n)(0)=c_i(\w)(0)\,,
\end{align*}
and as before this implies that $c_i(\w)=x_i$, as wanted.

Let $(\w,x)\in K$ and define
$x_\alpha=(1-\alpha)\,a(\w)+\alpha\,x\in B_{\widehat k_0}\subset BU$ for
$\alpha\in[0,1]$,~and
\[L=\{\alpha\in[0,1]\mid \mathcal{O}(\w,x_\alpha)=K^*\}\,.\]
If we prove that $L=[0,1]$, then $K=K^*$ and the proof is
finished. From the monotone character of the semiflow and since
$\mathcal{O}(\w,a(\w))=K^*$, it is immediate to check that if
$0<\alpha\in L$ then $[0,\alpha]\subset L$.

Next we show that $L$ is closed, that is, if $[0,\alpha)\subset L$
then $\alpha\in L$. Since $\{\tau(t,\w,x_\alpha)\mid t\geq 0\}$ is
uniformly stable for $\leq_A$ in bounded sets, let
$\delta(\varepsilon)>0$ be the modulus of uniform stability for
$\varepsilon>0$ in $B_{\widehat k_0}$. Thus, we take $\beta\in[0, \alpha)$
with $\di(x_\alpha,x_\beta)<\delta(\varepsilon)$  and we obtain
$\di(u(t,\w,x_\alpha),u(t,\w,x_\beta))<\varepsilon$ for each
$t\geq 0$. Moreover, $\mathcal{O}(\w,x_\beta)=K^*$ and hence,
there is a $t_0$ such that
$\di(u(t,\w,x_\beta),c(\w{\cdot}t))<\varepsilon$ for each $t\geq t_0$.
Then, we deduce that
$\di(u(t,\w,x_\alpha),c(\w{\cdot}t))<2\,\varepsilon$ for each $t\geq
t_0$ and $\mathcal{O}(\w,x_\alpha)=K^*$, i.e. $\alpha\in L$, as
claimed.

Finally, we prove that the case
$L=[0,\alpha]$ with $0\leq\alpha<1$ is impossible.  For each $i\in J$ we consider
the continuous maps
\begin{align*}
 K  \longrightarrow (0,\infty) &\,,\quad (\wt \w,\wt x)\mapsto
\wt x_i(0)-c_i(\wt \w)(0)\,,\\
K  \longrightarrow (0,\infty)& \,,\quad (\wt \w,\wt x)\mapsto
F_i(\wt \w,\wt x) -F_i(\wt \w, c(\wt \w))- (A(D\wt x- Dc(\wt \w)))_i\,.
\end{align*}
As explained above $\wt x_i(0)-c_i(\wt\w)(0)>0$. Moreover, from
$c(\wt\w)\leq_A \wt x$ and (F6) we deduce that $F_i(\wt
\w,\wt x) -F_i(\wt \w, c(\wt \w))-(A(D\wt x- Dc(\wt \w)))_i>0$. Hence,
there is an $\ep>0$ such that $\wt x_i(0)-c_i(\wt \w)(0)\geq \ep$ and
$F_i(\wt \w,\wt x) -F_i(\wt \w, c(\wt \w))-(A(D\wt x- Dc(\wt
\w)))_i>\ep$ for each $(\wt \w,\wt x)\in K$. Besides, since $(\wt
\w{\cdot}s,u(s,\wt \w,\wt x))\in K$, $u(s,\wt \w,\wt x)(0)=\wt x(s)$ for
each $s\leq 0$ because $K$ admits a flow extension, and $c(\wt
\w)(s)=c(\wt \w{\cdot}s)(0)$, we deduce that $\wt x_i(s)-c_i(\wt \w)(s)\geq
\ep$ for each $s\in (-\infty,0]$ and $(\wt \w,\wt x)\in K$.

It is not hard to check that $\cup_{\beta\in[0,1]}\cls_{\Om\times
X}\{\tau(t,\w,x_\beta)\mid t\geq 0\}$
 is a compact set. Hence,  since $\{\tau(t,\w,x_\alpha)\mid
t\geq 0\}$ is uniformly stable for $\leq_A$ in $B_{\widehat k_0}$,
$(\w,x_\beta)\in B_{\widehat k_0}$ for each $\beta\in[0,1]$,
 and hypotheses (D2) and (F3) hold, we deduce that there is
a $\delta>0$ such that
$\n{u(t,\w,x_\gamma)(0)-u(t,\w,x_\alpha)(0)}<\ep/4$ and
  \[\n{F(\w{\cdot}t,
u(t,\w,x_\gamma))-F(\w{\cdot}t, u(t,\w,x_\alpha))- A(D
u(t,\w,x_\gamma)- Du(t,\w,x_\alpha))}<\frac{\ep}{4}\] for each
$t\geq 0$ and $\gamma\in (\alpha,1]$ with
$\di(x_\alpha,x_\gamma)<\delta$. Besides, $\alpha\in L$, i.e.
$\mathcal{O}(\w,x_\alpha)=K^*$  and there is a $t_0\geq 0$ such
that $\n{u(t,\w,x_\alpha)(0)-c(\w{\cdot}t)(0)}<\ep/4$ and
  \[\n{F(\w{\cdot}t,
u(t,\w,x_\alpha))-F(\w{\cdot}t, c(\w{\cdot}t))- A(D u(t,\w,x_\alpha)- D
c(\w{\cdot}t))}<\frac{\ep}{4}\] for each $t\geq t_0$. Consequently, for
each $t\geq t_0$
\begin{equation}\label{xgamma}
\begin{array}{c}
\n{u(t,\w,x_\gamma)(0)-c(\w{\cdot}t)(0)}<{\displaystyle \frac{\ep}{2}}\,,\\[.2cm]
\n{F(\w{\cdot}t, u(t,\w,x_\gamma))-F(\w{\cdot}t, c(\w{\cdot}t))-A (D
u(t,\w,x_\gamma)- D c(\w{\cdot}t))}<{\displaystyle \frac{\ep}{2}}\,.
\end{array}
\end{equation}
Let $(\wt\w,\wt x)\in \mathcal{O}(\w,x_\gamma)$, i.e.~$(\wt\w,\wt
x)=\lim_{\,n\to\infty}(\w{\cdot}t_n,u(t_n,\w,x_\gamma))$ for some
$t_n\uparrow\infty$. The monotonicity and $c(\w)\leq_A x_\gamma$
imply that $c(\w{\cdot}t_n)\leq_A u(t_n,\w,x_\gamma)$, which yields
$c(\wt\w)\leq_A \wt x$.
From~\eqref{xgamma} there is an $n_0$
such that for each $i\in\{1,\ldots,m\}$
\[
0\leq u_i(t_n,\w,x_\gamma)(0)-c_i(\w{\cdot}t_n)(0)<\ep/2\,,\]
\[
F_i(\w{\cdot}t_n, u(t_n,\w,x_\gamma))-F_i(\w{\cdot}t_n, c(\w{\cdot}t_n))-(A(D
u(t_n,\w,x_\alpha)- D c(\w{\cdot}t_n)))_i<\ep/2
\]
for each $n\geq n_0$. Hence,
$0\leq \wt x_i(0)-c_i(\wt\w)(0)\leq\ep/2$ and
\[0\leq F_i(\wt
\w, \wt x)-F_i(\wt \w, c(\wt\w))-(A(D \wt x- D c(\wt \w)))_i\leq{\displaystyle \frac{\ep}{2}}\,.\]
As before, since this is true for each $(\wt\w,\wt
x)\in\mathcal{O}(\w,x_\gamma)$, admitting a flow extension, and $(\wt\w{\cdot}s,\wt x_s)\in \mathcal{O}(\w,x_\gamma)$ we
deduce that
\begin{equation}\label{epsilonmedio}
\begin{array}{c}
0\leq \wt x_i(s)-c_i(\wt\w)(s)\leq {\displaystyle \frac{\ep}{2}}\,,\\[.2cm]
0\leq F_i(\wt\w{\cdot}s,\wt x_s)-F_i(\wt\w{\cdot}s,c(\wt\w{\cdot}s))- (A(D\wt x_s-Dc(\wt\w{\cdot}s)))_i\leq {\displaystyle \frac{\ep}{2}}
\end{array}
\end{equation}
for each
$s\in(-\infty,0]$ and $i\in\{1,\ldots,m\}$.
Given any $(\wt\w,z)\in K$,  as shown above,
\begin{equation*}
\begin{array}{c}
z_i(s)-c_i(\wt\w)(s)\geq\ep\,,\\[.2cm]
F_i(\wt\w{\cdot}s,z_s)-F_i(\wt\w{\cdot}s,c(\wt\w{\cdot}s))-(A(Dz_s-Dc(\wt\w{\cdot}s)))_i\geq \ep
\end{array}
\end{equation*}
for each $s\in(-\infty,0]$ and $i\in J$, which combined with~\eqref{epsilonmedio} yields
\begin{equation}\label{inJ}
\wt x_i(s)\leq z_i(s)\quad \text{ and } \quad
F_i(\wt\w{\cdot}s,z_s)-F_i(\wt\w{\cdot}s,\wt x_s)-(A(Dz_s-D\wt x_s))_i\geq 0
\end{equation}
for each $s\in(-\infty,0]$ and $i\in J$. In addition, we deduce from Proposition~\ref{clasec1} that $c(\w)$, $z$ and $\wt x\in C^1(\R,\R^m)$.

Next, we claim that if $i\notin J$ the equality holds in~\eqref{inJ}. We know that $c_i(\wt\w)=z_i$,
\[F_i(\wt\w,z)-F_i(\wt\w,c(\wt\w))\geq (A(Dz-Dc(\wt\w)))_i\,,\]
and $(A(z-c(\wt\w)))_i=0$ because $0=z_i'-c_i(\wt\w)'\geq (A(z-c(\wt\w)))_i\geq 0$. Hence, from (D7) and~\eqref{Dformula} we deduce that
\begin{equation}\label{notinJ1}
(A(Dz-Dc(\wt\w)))_i=-\int_{-\infty}^0\left([d\nu(s)]\,A(z(s)-
c(\wt\w)(s))\right)_i\,.
\end{equation}
Moreover, $z'-c(\wt\w)'\geq A(z-c(\wt\w))$ and (D5) yield
\[ \int_{-\infty}^0 \left(\,[d\nu(s)](z'(s)-c(\wt\w)'(s))\right)_i\geq
\int_{-\infty}^0\left(\, [d\nu(s)]\,A(z(s)-c(\wt\w)(s))\right)_i\,,\]
which together with~\eqref{notinJ1} provide
\begin{multline*}
(A(Dz-Dc(\wt\w)))_i\geq - \int_{-\infty}^0 \left(\,[d\nu(s)](z'(s)-c(\wt\w)'(s))\right)_i=\frac{d}{ds} \left(D(z-c(\wt\w))\right)_i\\=F_i(\wt\w,z)-F_i(\wt\w,c(\wt\w))
\geq (A(Dz-Dc(\wt\w)))_i\,,
\end{multline*}
that is,
\[F_i(\wt\w,z)-F_i(\wt\w,c(\wt\w))
= (A(Dz-Dc(\wt\w)))_i\,.\]
From $c(\wt\w)\leq_A\wt x$ and $c_i(\wt\w)=\wt x_i$ we also deduce that
\[F_i(\wt\w,\wt x)-F_i(\wt\w,c(\wt\w))
= (A(D\wt x-Dc(\wt\w)))_i\,,\]
to conclude that $F_i(\wt\w,z)-F_i(\wt\w,\wt x)
= (A(Dz-D\wt x))_i$. Hence, if $i\notin J$ and $s\leq 0$
\[\wt x_i(s)=z_i(s)\quad \text{and} \quad F_i(\wt\w{\cdot}s,z_s)-F_i(\wt\w{\cdot}s,\wt x_s)
= (A(Dz_s-D\wt x_s))_i\,,\]
as stated, which together with~\eqref{inJ} provide $\wt x\leq z$ and
\[\frac{d}{ds} D z_s- \frac{d}{ds}D \wt x_s -A(Dz_s-D\wt x_s)\geq 0\,,\quad s\leq 0\,,\]
that is,  $\widehat D (z'-A z- (\wt x'-A\wt x))\geq 0$  from (D7). Finally,
since $\widehat D^{-1}$ is obviously positive in Cases I-II, and in Case
III because  of (D6), we conclude that $\wt x'-A
x\leq z'-A z$, that is, $\wt x\leq_A z$.

Since this holds for each $(\wt\w,z)\in K$, the definition of $a$
yields $c(\wt\w)\leq_A \wt x\leq_A a(\wt\w)$.
From~\eqref{copia} we know that $\mathcal{O}(\wt\w,a(\wt\w))=K^*$
and therefore $\mathcal{O}(\wt\w,\wt x)=K^*\subseteq
\mathcal{O}(\w, x_\gamma)$. Finally, from (F5) and
Theorem~\ref{minimal} we know that $\mathcal{O}(\w,x_\gamma)$ is a
minimal set, and we conclude that $\mathcal{O}(\wt\w,\wt
x)=\mathcal{O}(\w,x_\gamma)=K^*$, that is,  $\gamma\in L$, a
contradiction. Therefore, $L=[0,1]$ and
$\mathcal{O}(\w_0,x_0)=K^*$, as stated.
\qquad\end{proof}

\section{Applications to compartmental systems}\label{example}
In this Section we apply the previous results to the study of
compartmental models for the mathematical description
of processes in which the transport of material among
compartments takes a non-negligible length of time, and each
compartment produces or swallows material.

Let us suppose that we have a system formed
by $m$ compartments $C_1,\ldots,C_m$, and denote by $z_i(t)$ the amount of
material within compartment $C_i$ at time $t$  for each
$i\in\{1,\ldots,m\}$. Material flows from compartment $C_j$ into
compartment $C_i$ through a pipe having a transit time
distribution given by a positive regular Borel measure $\mu_{ij}$
with finite total variation  $\mu_{ij}(-\infty,0]=1$, for each
$i$, $j\in\{1,\ldots,m\}$,  whereas the
outcome of material from $C_i$ to $C_j$ is assumed to be instantaneous.   Let $\wt
g_{ij}:\mathbb{R}\times\mathbb{R}\to\mathbb{R}$ be the so-called
{\em transport function\/} determining the volume of material
flowing from $C_j$ to $C_i$ given in terms of the time $t$ and the
value of $z_j(t)$ for $i,j\in\{1,\ldots,m\}$. For each
$i\in\{1,\ldots,m\}$, at time $t\ge 0$, the compartment $C_i$
produces material itself at a rate $\int_{-\infty}^0
z'_i(t+s)\,d\nu_i(s)$, where $\nu_i$ is a positive regular
Borel measure with finite total variation
$\nu_i(-\infty,0]<\infty$ and $\nu_i(\{0\})=0$.
We will assume that the system is closed, that is, there is not inflow or outflow
of material from or to the environment surrounding the system.

 Once  the destruction and creation of
material is taken into account, the change of the amount of
material of any compartment $C_i$, $1\leq i\leq m$, equals the
difference between the amount of total influx into and total
outflux out of $C_i$, and we obtain a model governed by the
following system of  NFDEs:
\begin{multline}\label{ejemplo}
\frac{d}{dt}\left[z_i(t)-\int_{-\infty}^0
z_i(t+s)\,d\nu_i(s)\right]=-\sum_{j=1}^m
\wt g_{ji}(t,z_i(t)) \\ + \sum_{j=1}^m\int_{-\infty}^0\wt
g_{ij}(t+s,z_j(t+s)) \,d\mu_{ij}(s)\,,
\end{multline}
$i=1,\ldots,m$, where $\wt g=(\wt
g_{ij})_{i,j}:\R\times\R\to \R^{m\times m}$, and we will assume that

\vspace{.1cm}
\begin{itemize}
\item[(H1)] $\wt g$ is $C^1$-{\em admissible}, that is, $\wt g$ is $C^1$ in its second variable
and $\wt g$, $\frac{\partial}{\partial v}
 \wt g$ are uniformly continuous and bounded
on $\R\times \{v_0\}$ for all $v_0\in\R$;
\item[(H2)] $\wt g_{ij}(t,0)=0$  and  $\wt g_{ij}(t,v)$ is increasing in $v$ for each $t\in\R$ and $i,j=1,\ldots,m$;
\item[(H3)] $\wt g$ is a recurrent function, i.e. its {\em hull} is minimal;
\item[(H4)]$\nu_i(\{0\})=0$, $\sum_{i=1}^m \nu_i(-\infty,0]<1$,
$\mu_{ij}(-\infty,0]=1$ and $\int_{-\infty}^0
|s|\,d\mu_{ij}(s)<\infty$ for $i,j=1,\ldots,m$;
\item[(H5)] for each $i=1,\ldots,m$,
 $L_i^+=\sum_{j=1}^m l_{ji}^+<+\infty$ and there is  a $\beta_i>0$ such that $ \beta_i\left(1-\int_{-\infty}^0e^{-\beta_i s}\,d\nu_i(s)\right)>L_i^+$
where $\;\;{\displaystyle l_{ji}^+=\sup_{{(t,\,v)\in\R^2}}\frac{\partial \wt g_{ji}}{\partial v}(t,v)}\;$.
\end{itemize}
As usual, we include the nonautonomous equation~\eqref{ejemplo} into a family of nonautonomous NFDEs as follows.
Let $\Om$ be the {\em hull\/} of $\wt g$, namely, the closure of
the set of mappings $\{\wt g_t\mid t\in \R\}$, with $\wt
g_t(s,v)=\wt g(t+s,v)$, $(s,v)\in \R^2$, with the topology of
uniform convergence on compact sets, which from (H1) is a compact
metric space (see Hino, Murakami and Naito~\cite{hino}).  Let $(\Om,\sigma,\R)$ be the
continuous flow~defined on $\Om$ by translation, $\sigma:\R\times
\Om\to \Om$, $(t,\w)\mapsto\w{\cdot}t$, with $\w{\cdot}t(s,v)= \w(t+s,v)$. From
hypothesis (H3), the flow $(\Om,\sigma,\R)$ is minimal. In
addition, if $\wt g$ is almost periodic (resp. almost automorphic)
the flow will be almost periodic (resp. almost automorphic).
Notice that these two cases are included in our formulation.

Let $g:\Om\times\R\to\R^{m\times m},\;
(\w,v)\mapsto\w(0,v)$, continuous on $\Om\times \R$ and denote
$g=(g_{ij})_{i,j}$. Let $F:\Om\times BU\to\R^m$ be
the map defined by
\begin{equation}\label{defiF}
F_i(\w,x)=-\sum_{j=1}^m
g_{ji}(\w,x_i(0))+\sum_{j=1}^m\int_{-\infty}^0g_{ij}(\w{\cdot}s
,x_j(s))\, d\mu_{ij}(s)\,,
\end{equation}
for $(\w,x)\in\Om\times BU$ and $i\in\{1,\ldots,m\}$. Hence, the
family
\addtocounter{equation}{1}
\begin{equation}\tag*{(\arabic{section}.\arabic{equation})$_\w$}\label{hejemplo}
 \frac{d}{dt}Dz_t=F(\w{\cdot}t,z_t)\,, \quad t\geq 0\,,\;\w\in\Om\,,
\end{equation}
is of the form~\ref{Ninfdelay} with $D_ix=x_i(0)-\int_{-\infty}^0
x_i(s)\,d\nu_i(s)$, $i=1,\ldots,m$
and includes system~\eqref{ejemplo} when $\w=\wt g$.

In~\cite{MNO2008}, Mu\~{n}oz-Villarragut, Novo and Obaya study NFDEs~\ref{hejemplo} which are monotone for the standard ordering.
This requires the positivity of the measures $d\eta_{i}=l_{ii}^-\,d\mu_{ii}-\sum_{k=0}^m
l_{ki}^+\,d\nu_{i}$, where
 \[l_{ii}^-=\inf_{(t,\,v)\in\R^2}\frac{\partial \wt
g_{ii}}{\partial v}(t,v)\,,\; \text{ and }\;
l_{ij}^+=\sup_{{(t,\,v)\in\R^2}}\frac{\partial \wt g_{ij}}{\partial
v}(t,v)\,,\quad i,j=1,\ldots m\,.\]
We now remove this condition, which controls the material
produced in the compartments in terms of the material transported through the pipes.
In order to verify the monotonicity for the exponential ordering, we assume hypothesis (H5),
which imposes a strong restriction on the size of the neutral term.
The next statement shows that the main conclusions obtained in~\cite{MNO2008} remain
valid for~\ref{hejemplo}.

\begin{teor}\label{teorejemplo1}
  Assume that {\rm (H1)-(H5)} hold and let $(\w_0,x_0)\in \Om\times BU$ with $x_0$ Lipschitz.
  Then the solution $z(t,\w_0,x_0)$ of~{\rm \ref{hejemplo}}$\!_{_0}$ with initial value $x_0$ is bounded and its omega-limit set is a copy of the base, that is, $\mathcal{O}(\w_0,x_0)=\{(\w,c(\w))\mid \w\in\Om\}$ for a continuous equilibrium
  $c\colon \Om\to BU$ and
  {\rm \[\lim_{t\to\infty}
\di(u(t,\w_0,x_0),c(\w_0{\cdot}t))=0\,.\]}
\end{teor}

\begin{proof}
It is easy to check that assumptions (D1)-(D6), and (F1)-(F3) hold, and we consider the local skew-product semiflow~\eqref{Ndelaskewcoop} induced by~\ref{hejemplo}.
Next we check the monotonicity assumptions (F4) and (F6) for the exponential ordering $\leq_{A,\infty}$ where $A$ is the quasipositive diagonal matrix with diagonal
entries $(-\beta_1,\ldots,-\beta_m)$
 given in (H5). Notice that (A1) and (D7) also hold. As usual $\leq_{A,\infty}$ will be denoted~by~$\leq_A$.

Let $x$, $y\in BU$ with $x\leq_A y$, that is,
$x\le y$ and $ y(t)-x(t)\geq e^{A(t-s)}(y(s)-x(s))$ if $-\infty<s\leq t\leq 0$. Then, from (H2) $g_{ij}$ are increasing in their second variable, and from~\eqref{defiF} we deduce that
\begin{multline*}
F_i(\w,y)-F_i(\w,x)+ \beta_i(D_iy- D_ix)\geq  \sum_{j=1}^m [g_{ji}(\w,x_i(0))-g_{ji}(\w,y_i(0))]\\ +\beta_i \left(y_i(0)-\int_{-\infty}^0
y_i(s)\,d\nu_i(s)-x_i(0)+ \int_{-\infty}^0
x_i(s)\,d\nu_i(s)\right)
\end{multline*}
for each $i=1,\ldots,m$.
In addition, $g_{ji}(\w,x_i(0))-g_{ji}(\w,y_i(0))\geq -l_{ji}^+\,(y_i(0)-x_i(0))$ and
 from $x\leq_A y$ we deduce that for each $s\leq 0$
 \[(y_i(s)-x_i(s))\,e^{\beta_i s}\leq y_i(0)-x_i(0)\,,\]
which yields
\begin{multline}\label{desih51}
F_i(\w,y)-F_i(\w,x)+\beta_i (D_iy- D_ix) \\ \geq \left[-L_i^++ \beta_i \left(1-\int_{-\infty}^0\,e^{-\beta_i s}\,d\nu_i(s)\right)\right](y_i(0)-x_i(0))\,,
\end{multline}
for each $i=1,\ldots,m$. Hence from hypothesis (H5) we deduce that (F4) and (F6) hold, as claimed.

In order to check hypothesis (F5) we define  $M\colon \Om\times BU\to \R$, the {\em total mass\/}  of the system~{\rm
\ref{hejemplo}} as
\[
M(\w,x):=\sum_{i=1}^m D_ix+  \sum_{i=1}^m \sum_{j=1}^m
\int_{-\infty}^0 \left( \int_s^0g_{ji}(\w{\cdot} \tau,x_i(\tau))\,d\tau
\right)d\mu_{ji}(s)\,,
\]
which is well defined from condition (H4).
It is easy to check, as shown in~\cite{MNO2008}, that $M$ is a uniformly continuous function on all the sets of the form $\Om\times B_k$ with $k>0$ for the product metric topology. Moreover, it is constant along the trajectories because from~\ref{hejemplo}
\[\frac{d}{dt}M(\tau(t,\w,x))=\frac{d}{dt}M(\w{\cdot}t,z_t(\w,x))=0\,,\]
that is, $M(\w{\cdot}t,z_t(\w,x))=M(\w,x)$  for each $t\geq 0$ where $z(t,\w,x)$ is defined.

Let $x$, $y\in BU$ with $x\leq_A y$. We can apply Theorem~\ref{monotone} to deduce that the induced semiflow is monotone and hence, $u(t,\w,x)\leq_A u(t,\w,y)$
whenever they are  defined. Thus, since the transport functions $g_{ji}$ are monotone and the measures $\mu_{ji}$ are positive
\[\int_{-\infty}^0 \!\! \left( \int_s^0\!\!\! (g_{ji}(\w{\cdot}(t+ \tau),z_i(t+\tau,\w,y)-g_{ji}(\w{\cdot} (t+\tau),z_i(t+\tau,\w,x))\,d\tau
\!\right)\!d\mu_{ji}(s)\geq 0\]
and we deduce that
\begin{align*}
\sum_{i=1}^m\left (D_iz_t(\w,y)-D_iz_t(\w,x)\right)& \leq M(\w{\cdot}t,z_t(\w,y))-M(\w{\cdot}t,z_t(\w,x))\\ & =M(\w,y)-M(\w,x)\,.
\end{align*}
Next we check that given $\varepsilon >0$ there is a $\delta>0$  such that if $x$, $y\in B_k$ with $x\leq_A y$ and
 $\di(x,y)<\delta$ then  $\n{z(t,\w,y)-z(t,\w,x)}\leq \varepsilon$ whenever they
are  defined.
From hypothesis (H4) we define $\gamma=\sum_{i=1}^m \nu_i(-\infty,0]<1$. From the uniform continuity of $M$, given
$\varepsilon_0=\varepsilon\,(1-\gamma)>0$ there exists  $0<\delta <\varepsilon_0$ such that if
 $x$, $y\in B_k$ with $\di(x,y)<\delta$ then
 $|M(\w,y)-M(\w,x)|<\varepsilon_0$. Therefore, if $x,\,y\in B_k$ and $x\leq_A y$
we have $\sum_{i=1}^m (D_iz_t(\w,y)-D_iz_t(\w,x))<\varepsilon_0$. Consequently, from the definition of $D_i$,
\begin{align*}
0\leq z_i(t,\w,y)-z_i(t,\w,x) & \leq \sum_{j=1}^m(z_j(t,\w,y)-z_j(t,\w,x))\\ & <\varepsilon_0+ \sum_{j=1}^m \int_{-\infty}^0 (z_j(t+s,\w,y)-z_j(t+s,\w,x))\,d\nu_j(s)\\
&\leq \varepsilon_0+ \gamma \,\n{z_t(\w,y)-z_t(\w,x)}_\infty  \,,
\end{align*}
from which we deduce that
 $\n{z_t(\w,y)-z_t(\w,x)}_\infty
(1-\gamma)<\varepsilon_0=\varepsilon(1-\gamma)$, that is,  $\n{z(t,\w,y)-z(t,\w,x)}\leq \varepsilon$ whenever they
are  defined, as claimed.

 From this and since 0 is a solution of~\ref{hejemplo} for each $\w\in\Om$  because $g_{ij}(\w,0)=0$ for each $i,j=1,\ldots,m$, we deduce that each solution $z(t,\w,x)$ with $x\geq_A 0$ or $x\leq_A 0$ is globally defined and bounded. As a consequence, and together with the monotone character of the semiflow, we conclude that all the solutions with Lipschitz initial data are globally defined and bounded because given $x\in BU$ and Lipschitz, it is easy to find $c_1<0$ and $c_2>0$ such that $c_1\,{\bf 1}\leq_A x\leq_A c_2\,{\bf 1}$.

 Let $(\w,x)\in\Om\times BU$ with $x$ Lipschitz, and $k>0$ such that $z_t(\w,x)\in B_{k}$
for all $t\geq 0$. Then, as above,
we deduce that given $\varepsilon>0$ there exists a $\delta>0$
such that
\[\n{z(t+s,\w,x)-z(t,\w{\cdot}s,y)}=\n{z(t,\w{\cdot}s,z_s(\w,x))-z(t,\w{\cdot}
s,y)}<\varepsilon
\]
for all $t\geq 0$ whenever $y\in B_k$, $y\leq_A z_s(\w,x)$ or $z_s(\w,x)\leq_A y$ and {\rm
$\di(z_s(\w,x),y)<\delta$}, which shows the uniform stability of
the trajectories for the order $\leq_A$ in $B_k$ for each $k>0$ and hypothesis (F5) holds.
Hence, the result follows from Theorem~\ref{maincopy}.
\qquad\end{proof}

Concerning the solutions of the original compartmental system, we
obtain the following result. Although the theorem is stated in the almost periodic case,
similar conclusions are obtained changing almost periodicity for
periodicity, almost automorphy or recurrence, that is, all
solutions are asymptotically of the same type as the transport
functions.

\begin{teor}
Under Assumptions {\rm (H1-H5)} and in the almost periodic case,
there are infinitely many almost periodic solutions of system~{\rm\eqref{ejemplo}} and
all the solutions with Lipschitz initial data are asymptotically almost periodic.
\end{teor}

\begin{proof}
Let $\w_0=\wt g$ and a Lipschitz initial value $x_0$. The omega-limit set  of each solution
$z(t,\w_0,x_0)$ is a copy of the base
$\mathcal{O}(\w_0,x_0)=\{(\w,x(\w))\mid \w\in\Om\}$ and hence,
$z(t,\w_0,x(\w_0))=x(\w_0{\cdot}t)(0)$ is an almost periodic solution
of~\eqref{ejemplo} with \[
\lim_{t\to\infty}\n{z(t,\w_0,x_0)-z(t,\w_0,x(\w_0))}=0\,.\]
Next we check that there are infinitely many minimal subsets. From the definition of the mass and (H4) \[M(\w_0,k\,{\bf 1})=\sum_{i=1}^m k\,(1-\nu_i(-\infty,0])+\sum_{i=1}^m \sum_{j=1}^m
\int_{-\infty}^0 \left( \int_s^0g_{ji}(\w_0{\cdot} \tau,k)\,d\tau
\right)d\mu_{ji}(s)\]
tends to $+\infty$ as $k$ goes to $+\infty$, and thus, given $c>0$ there is a
$k_c>0$ such that $M(\w_0,k_c\,{\bf 1})=c$. Finally,  since the mass is constant along the trajectories, $\mathcal{O}(\w_0,k_c{\bf 1})$ provides a different minimal set and hence a different almost periodic solution for each $c>0$.
\qquad\end{proof}

Next we consider the particular case where the compartmental systems are described by NFDEs with finite delay
\begin{equation}\label{ejemplo2}
\frac{d}{dt}(z_i(t)-\gamma_i\,z_i(t-\alpha_i))=-\sum_{j=1}^m
\wt g_{ji}(t,z_i(t))+ \sum_{j=1}^m
\wt g_{ij}(t-\rho_{ij},z_j(t-\rho_{ij}))\,,
\end{equation}
$i=1,\ldots,m$, where $\wt g=(\wt
g_{ij})_{i,j}:\R\times\R\to \R^{m\times m}$, and we will assume that

\vspace{.1cm}
\begin{itemize}
\item[(G1)] $\wt g$ is $C^1$-{\em admissible}, that is, $\wt g$ is $C^1$ in its second variable
and $\wt g$, $\frac{\partial}{\partial v}
 \wt g$ are uniformly continuous and bounded
on $\R\times \{v_0\}$ for all $v_0\in\R$;
\item[(G2)] $\wt g_{ij}(t,0)=0$  and  $\wt g_{ij}(t,v)$ is increasing in $v$ for each $t\in\R$ and $i,j=1,\ldots,m$;
\item[(G3)] $\wt g$ is a recurrent function, i.e. its {\em hull} is minimal;
\item[(G4)] $\sum_{i=1}^m \gamma_i<1$, $\alpha_i> 0$ and $\rho_{ij}\geq 0$, $i,j=1,\ldots,m$;
\item[(G5)] for each $i=1,\ldots,m$ one of the following hypotheses holds:
 \begin{itemize}
 \item[(G5.1)] $L_i^+=\sum_{j=1}^m l_{ji}^+<+\infty$ and there is  a $\beta_i>0$ such that \[ \beta_i(1-\gamma_i\,e^{\beta_i\alpha_i})>L_i^+\,;\]
 \item[(G5.2)]  $L_i^+=\sum_{j=1}^m l_{ji}^+<+\infty$, $\alpha_i\geq \rho_{ii}$  and there is a $\beta_i\geq L_i^+$ such that
 \[\beta_i\,(1-\gamma_i\,e^{\beta_i\alpha_i})+e^{\beta_i\rho_{ii}}\,l_{ii}^->L_i^+ \,,\]
 \end{itemize}
where $\;\;{\displaystyle l_{ij}^+=\sup_{{(t,\,v)\in\R^2}}\frac{\partial \wt g_{ij}}{\partial v}(t,v)}\;$  and  $\;{\displaystyle l_{ii}^-=\inf_{(t,v)\in \R^2}\frac{\partial\wt g_{ii}}{\partial v}(t,v)}$.
\end{itemize}

As before, the family
\addtocounter{equation}{1}
\begin{equation}\tag*{(\arabic{section}.\arabic{equation})$_\w$}\label{hejemplo2}
 \frac{d}{dt}Dz_t=F(\w{\cdot}t,z_t)\,, \quad t\geq 0\,,\;\w\in\Om\,,
\end{equation}
with
\begin{equation}\label{defiF2}
F_i(\w,x)=-\sum_{j=1}^m
g_{ji}(\w,x_i(0))+\sum_{j=1}^mg_{ij}(\w{\cdot}(-\rho_{ij})
,x_j(t-\rho_{ij}))\,,
\end{equation}
is of the form~\ref{Ninfdelay} with $D_ix=x_i(0)-\gamma_i \,x(-\alpha_i)$, $i=1,\ldots,m$
and includes system~\eqref{ejemplo} when $\w=\wt g$.

Under these assumptions we deduce the convergence of the solutions with Lipschitz initial data. Notice that hypothesis (G5.2), in this case of finite delay, improves the conditions of applicability of the theory
stated in Theorem~\ref{teorejemplo1}.

\begin{teor}
  Assume that {\rm (G1)-(G5)} hold and let $(\w_0,x_0)\in \Om\times BU$ with $x_0$ Lipschitz.
  Then the solution $z(t,\w_0,x_0)$ of~{\rm \ref{hejemplo2}}$\!_{_0}$ with initial value $x_0$ is bounded and its omega-limit set is a copy of the base, that is, $\mathcal{O}(\w_0,x_0)=\{(\w,c(\w))\mid \w\in\Om\}$ for a continuous equilibrium
  $c\colon \Om\to BU$ and
  {\rm \[\lim_{t\to\infty}
\di(u(t,\w_0,x_0),c(\w_0{\cdot}t))=0\,.\]}
\end{teor}

\begin{proof}
The only difference with the proof of Theorem~\ref{teorejemplo1} is to check (F4) and (F6)
for the components $i\in\{1,\ldots,m\}$ for which (G5.2) holds.

Thus, let $i\in\{1,\ldots,m\}$ satisfy (G5.2), and let $x$, $y\in BU$ with $x\leq_A y$.
Then, from~\eqref{defiF2}, $g_{ii}(\w,y(-\rho_{ii}))-g_{ii}(\w,x(-\rho_{ii}))\geq l_{ii}^-\,(y_i(-\rho_{ii})-x_i(-\rho_{ii}))$
and \[(y_i(-\alpha_i)-x_i(-\alpha_i))\,e^{-\beta_i\alpha_i}\leq (y(-\rho_{ii})-x(-\rho_{ii}))e^{-\beta_i\rho_{ii}}\,,\]
because $\alpha_i\geq \rho_{ii}$ and $x\leq_Ay$, we deduce that
\begin{multline*}
F_i(\w,y)-F_i(\w,x)+\beta_i (D_iy- D_ix)
\geq  \,(\beta_i-L_i^+)\,(y_i(0)-x_i(0))\\ +(l_{ii}^--\beta_i\,\gamma_i\,e^{\beta_i(\alpha_i-\rho_{ii})})\,(y_i(-\rho_{ii})-
x_i(-\rho_{ii}))\,.
\end{multline*}
In addition, from $\beta_i\geq L_i^+$, $(y(-\rho_{ii})-x(-\rho_{ii}))\,e^{-\beta_i\rho_{ii}}\leq y_i(0)-x_i(0)$
and (G5.2) we conclude that
\begin{multline*}
F_i(\w,y)-F_i(\w,x)+\beta_i (D_iy- D_ix)\geq (y_i(-\rho_{ii})-x_i(-\rho_{ii}))\big[(\beta_i-L_i^+)\,e^{-\beta_i\rho_{ii}}\\ +l_{ii}^--\beta_i\,\gamma_i\,e^{\beta_i(\alpha_i-\rho_{ii})}\big]\,,
\end{multline*}
and (F4) and (F6) hold.
\qquad\end{proof}

As before, concerning the solutions of the original compartmental system, we
obtain the following result.

\begin{teor}
Under Assumptions {\rm (G1-G5)} and in the almost periodic case,
there are infinitely many almost periodic solutions of system~{\rm\eqref{ejemplo2}} and
all the solutions with Lipschitz initial data are asymptotically almost periodic.
\end{teor}


\begin{thebibliography}{99}
\bibitem{ARI1987} {\sc O. Arino, E. Haourigui}, On the asymptotic
        behavior of solutions of some delay differential systems which
        have  a first integral, {\em J. Math. Anal. Appl.}
        \textbf{122} (1987), 36--46.
\bibitem{AB1990} {\sc O. Arino, F. Bourad}, On the asymptotic
        behavior of the solutions of a class of scalar neutral
        equations generating a monotone semiflow, {\em J. Differential
        Equations} \textbf{87} (1990), 84--95.
\bibitem{aufr} {\sc J.-P. Aubin, H. Frankowska},
        {\em Set-Valued Analysis},
        Birkh\"{a}user, Boston, Basel, Berlin, 1990.
\bibitem{choq} {\sc G. Choquet},
        {\em Lectures on Analysis. Integration and Topological Vector Spaces,
         Vol. I\/}, Math. Lecture Notes, Benjamin, Reading, MA, 1969.
\bibitem{chue} {\sc I.D. Chueshov}, \textit{Monotone Random
        Systems. Theory and Applications}, Lecture Notes in Math. \textbf{1779},
        Springer-Verlag, Berlin, Heidelberg, 2002.
\bibitem{elli} {\sc R. Ellis}, \textit{Lectures on Topological
        Dynamics\/}, Benjamin, New York, 1969.
\bibitem{gyell} {\sc I. Gy\"{o}ri, J. Eller}, Compartmental systems with pipes,
        {\em Math. Biosci.}  {\bf 53} (1981) 223--247.
\bibitem{gyori} {\sc I. Gy\"{o}ri}, Connections between compartmental
        systems with pipes and integro-differential equations, {\em  Math.
        Modelling\/}  {\bf 7}  (1986) 1215--1238.
\bibitem{gyoriwu} {\sc I. Gy\"{o}ri, J. Wu}, A neutral equation arising
        from compartmental systems with pipes, {\em J. Dynamics
        Differential Equations} \textbf{3} No.2
        (1991), 289--311.
\bibitem{hale} {\sc J.K. Hale},
        \textit{Theory of Functional Differential
        Equations}, Applied Mathematical Sciences {\bf 3},
        Springer-Verlag, Berlin, Heidelberg, New York 1977.
\bibitem{hale2} {\sc J.K. Hale, S.M. Verduyn Lunel},
        \textit{Introduction to Functional Differential
        Equations}, Applied Mathematical Sciences {\bf 99},
        Springer-Verlag, Berlin, Heidelberg, New York 1993.
\bibitem{hino} {\sc Y. Hino, S. Murakami, T. Naito},
        \textit{Functional Differential Equations with Infinite
        Delay}, Lecture Notes in Math. \textbf{1473},
        Springer-Verlag, Berlin, Heidelberg, 1991.
\bibitem{H1988} {\sc M.W. Hirsch}, Stability and convergence in strongly monotone dynamical
         systems, {\em J. Reine Angew. Math.\/} \textbf{383} (1988), 1--53.
\bibitem{JJ1996} {\sc J.A. Jacquez}, \textit{Compartmental
        Analysis in Biology and Medicine}, Third Edition, Thomson-Shore Inc., Ann
        Arbor, Michigan, 1996.
\bibitem{JS1993} {\sc J.A. Jacquez, C.P. Simon}, Qualitative theory of compartmental
        systems, {\em SIAM Review}  \textbf{35} No.1 (1993), 43--79.
\bibitem{JS2002} {\sc J.A. Jacquez, C.P. Simon}, Qualitative theory of compartmental systems
        with lags, {\em Math. Biosci.\/} \textbf{180} (2002), 329--362.
\bibitem{jizh} {\sc J. Jiang, X.-Q. Zhao}, Convergence in monotone
        and uniformly stable skew-product semiflows with
        applications, {\em J. Reine Angew. Math} \textbf{589}
        (2005), 21--55.
\bibitem{KR1992} {\sc T. Krisztin}, An invariance principle of Lyapunov-Razumikhin type and compartmental systems,  World Congress of Nonlinear Analysts' 92, Vol. I--IV (Tampa, FL, 1992),  1371--1379, de Gruyter, Berlin, 1996.
\bibitem{KWU1996} {\sc T. Krisztin, J. Wu}, Asymptotic Periodicity,
        Monotonicity, and Oscillation of Solutions of Scalar
        Neutral Functional Differential Equations, {\em J. Math.
        Anal. Appl.\/} \textbf{199} (1996), 502--525.
\bibitem{MA1984} {\sc H. Matano}, Existence of nontrivial unstable
        sets for equilibriums of strongly order preserving systems, {\em
        J. Fac. Sci. Univ. Kyoto\/} \textbf{30} (1984), 645--673.
\bibitem{MNO2008} {\sc V. Mu\~{n}oz-Villarragut, S. Novo, R. Obaya}, Neutral functional
        differential equations with applications to compartmental systems,
        \textit{SIAM J. Math. Anal.} \textbf{40} No. 3 (2008) 1003--1028.
\bibitem{NNO2005} {\sc S. Novo, C. N\'{u}\~{n}ez, R. Obaya}, Almost automorphic and
        almost periodic dynamics for quasimonotone
        non-autonomous functional differential equations,
        {\em J. Dynamics Differential Equations} \textbf{17} No.3
        (2005), 589--619.
\bibitem{NOS2007} {\sc S. Novo, R. Obaya, A.M. Sanz}, Stability and extensibility results
        abstract skew-product semiflows, {\em J. Differential Equations\/} \textbf{235} No. (2007), 623-646.
\bibitem{PO1989} {\sc P. Pol\'{a}\v{c}ik}, Convergence in smooth strongly monotone flows defined
        by semilinear parabolic equations,  {\em J.
        Differential Equations\/} {\bf 79} (1989), 89--110.
\bibitem{sase1} {\sc R.J. Sacker, G.R. Sell},
        \textit{Lifting properties in skew-products flows with
        applications to differential equations}, Mem. Amer. Math.
           Soc. \textbf{190}, Amer. Math. Soc., Providence 1977.
\bibitem{shyi} {\sc W. Shen, Y. Yi}, Almost Automorphic and Almost
        Periodic Dynamics in Skew-Product Semiflows, {\em Mem. Amer. Math. Soc.}
        {\bf 647}, Amer. Math. Soc., Providence  1998.
\bibitem{smith1995} {\sc H.L. Smith}, \textit{Monotone Dynamical Systems. An
        introduction to the Theory of Competitive and Cooperative Systems}, Amer. Math. Soc., Providence, 1995.
\bibitem{SMTH1990} {\sc H.L. Smith, H.R. Thieme}, Monotone semiflows in scalar non-quasi
        monotone functional differential equations, {\em J. Math. Anal. Appl.\/} \textbf{150} (1990), 289--306.
\bibitem{SMTH1991} {\sc H.L. Smith, H.R. Thieme}, Strongly order preserving semiflows generated by functional differential equations, {\em J. Differential Equations\/}
        \textbf{93} (1991), 332-363.
\bibitem{WW1985} {\sc Z. Wang, J. Wu}, Neutral Functional
        Differential Equations with Infinite delay, {\em Funkcial.
        Ekvac.\/} \textbf{28} (1985), 157--170.
\bibitem{jwu1991} {\sc J. Wu}, Unified treatment of local theory
        of NFDEs with infinite delay, {\em Tamkang J. Math.\/}
        \textbf{22} No. 1 (1991), 51--72.
\bibitem{WU1991} {\sc J. Wu}, Asymptotic periodicity of solutions to a class of neutral
        functional differential equations, {\em Proc. Amer. Math. Soc\/} \textbf{113} (1991), 355-363.
\bibitem{WF1991} {\sc J. Wu, H.I. Freedman}, Monotone semiflows generated by neutral
        functional differential equations with application to compartmental systems,
        {\em  Can. J. Math.\/} \textbf{43} (5) (1991), 1098-1120.
\bibitem{WZ2002} {\sc J. Wu, X.-Q. Zhao}, Diffusive monotonicity and threshold dynamics of delayed reaction diffusion equations {\em J. Differential Equations\/} \textbf{186} (2002), 470--484.
\end{thebibliography}
\end{document}